\documentclass[12pt]{article}
\usepackage{a4,latexsym,amssymb,amsmath}
\usepackage[all]{xy}
\usepackage[active]{srcltx}
\usepackage{tikz,graphicx}
\usetikzlibrary{arrows,automata,backgrounds,calendar}
\DeclareMathAlphabet{\mathpzc}{OT1}{pzc}{m}{it}

\newtheorem{theorem}{Theorem}[section]
\newtheorem{proposition}[theorem]{Proposition}
\newtheorem{corollary}[theorem]{Corollary}
\newtheorem{lemma}[theorem]{Lemma}
\newtheorem{definition}[theorem]{Definition}

\newcommand{\Z}{{\mathbb Z}}

\def\varprojlim{\mathop{\lim\limits_{{\longleftarrow}}\,}}
\def\itvarprojlim{\underset{\longleftarrow}{\it lim}}
\def\itvarprojlim#1{\underset{#1}{\underset{\longleftarrow}{\it lim}\,}}
\def\c{\mathfrak c}
\def\a{\mathfrak a}
\def\p{\mathfrak p}
\def\P{\mathfrak P}
\def\ds{\displaystyle}
\def\hang{\hangindent\Itemindent}
\def\textindent#1{\hskip\Itemindent\llap{\hbox{\rm #1}\enspace}\ignorespaces}
\def\Item{\par\noindent\hang\textindent}
        \newdimen\Itemindent \Itemindent=.9cm
\def\demo #1{\vskip-\lastskip \vskip0.8truecm\noindent{\bf #1}\enspace}
\def\noproof{{\unskip\nobreak\hfill\penalty50\hskip2em\hbox{}\nobreak\hfill$\square$\parfillskip=0pt\finalhyphendemerits=0\par}}
\def\enddemo{\ifmmode\eqno\square\else\noproof\vskip0.8truecm\fi}
\def\Cal#1{{\cal #1}}

\def\res{\hbox{\it res}}

\def\ttcup{\mathop{\mathord{\bigcup}\mkern-9.7mu^{\hbox{.}}\mkern 9mu}}
\def\ktcup{\mathop{\mathord{\cup}\mkern-8.5mu^{\hbox{.}}\mkern 9mu}}
\def\tcup{\mathop{\mathord{\bigcup}\mkern-12.7mu^{\hbox{.}}\mkern 9mu}}
\def\Tcup{\displaystyle \tcup}

\def\pfeil{\mathop{\rightarrow}\limits}

\def\lpfeil{\mathop{\longrightarrow}\limits}
\def\liso{\mathrel{\hbox{$\longrightarrow$} \kern-14pt\lower-4pt\hbox{$\scriptstyle\sim$}\kern7pt}}
\def\kiso{\mathrel{\hbox{$\rightarrow$} \kern-10.5pt\lower-4pt\hbox{$\scriptstyle\sim$}\kern4pt}}
\def\lr@iso{\mathrel{\kern6pt\lower2pt\hbox{$\scriptstyle\sim$} \kern-12pt\hbox{$\longrightarrow$}}}
\def\lriso#1{\mathrel{\mathop{\lr@iso}\limits^{#1}}}
\def\kr@iso{\mathrel{\kern6pt\lower2pt\hbox{$\scriptstyle\sim$} \kern-9pt\hbox{$\rightarrow$}}}
\def\kriso#1{\mathrel{\mathop{\kr@iso}\limits^{#1}}}
\def\surjr#1{\stackrel{#1}{\hbox{$\relbar \! \! \! \twoheadrightarrow$}}}
\def\injr#1{\stackrel{#1}{\hbox{$ \hookrightarrow$}}}

\def\Coind{{\mathop{\it Coind}}}
\font\emar = cmsy10 scaled\magstep4
\def\freeproductsmall{\mathop{\hbox{$\ast$}}\limits}
\def\freeproduct{\mathop{\lower.4mm\hbox{\emar \symbol{3}}}\limits}
\def\ffreeproduct{\mathop{\lower.4mm\hbox{\emar \symbol{3}}}}

\newcommand{\xyalign}{\entrymodifiers={+!!<0pt,\fontdimen22\textfont2>}}

\makeatletter
\newcommand{\superimpose}[2]{
  {\ooalign{$#1\@firstoftwo#2$\cr\hfil$#1\@secondoftwo#2$\hfil\cr}}}
\makeatother

\newcommand{\resprod}{\mathpalette\superimpose{{\prod}{\coprod}}}
\def\Resprod#1{\underset{#1}{\mathpalette\superimpose{{\hbox{$\prod$}}{\hbox{$\coprod$}}}}}
\def\sprod#1{\underset{#1}{\hbox{$\prod$}}}

\begin{document}

\title{\bf\boldmath  Corestricted Free Products\\ of  Profinite Groups}
\author{by Jochen G\"artner and  Kay Wingberg}
\date{
}
\maketitle

\begin{abstract}
\vspace{.1cm}
 \scriptsize We introduce the notion of corestricted free products of a family of profinite groups indexed over an arbitrary profinite space. 
 Using arithmetic results of the second author, this enables us to prove an analogue of Riemann's existence theorem for the decomposition groups of certain infinite sets of primes of a number field.
\end{abstract}
\normalsize

\indent

In 1971, Neukirch \cite{N} introduced the concept of free products of families of pro-$\c$-groups indexed over a discrete set. Since in number theory projective limits of free pro-$\c$-products arise in a natural way, it became necessary to generalize this concept to families varying continuously over a profinite index space. This has been worked out by several authors by introducing the notion of compact bundles $\Cal G$ of pro-$\c$-groups over a profinite space $T$, see \cite{GR}, \cite{Har}, \cite{Mel} and \cite{NSW}. These bundles are group objects in the category of profinite spaces over $T$ such that every fibre is a pro-$\c$-group. The generalized free products obtained in this way turned out to be the key tool to prove number theoretical analogues of the Riemann existence theorem over large number fields, i.e.\ inertia groups, indexed by  a projective limit of possibly infinite sets of primes of a tower of number fields, form a free product. Furthermore, analogous results hold for the decomposition groups lying over a {\it finite} set of primes of a number field, see \cite{NSW} chap.10 \S5.

From the arithmetic point of view, it is interesting whether infinitely many decomposition groups also form a free product. However the notion of compact bundles turns out to be too restrictive. 
Inspired by the abelian case in number theory, i.e.\ the idele group of a number field which is a restricted product and topological not compact, Neukirch introduced the concept of corestricted free pro- \linebreak ducts of pro-$\c$-groups over discrete sets. Generalizing this concept to a more general profinite space $T$ leads to possibly non-compact bundles, which are group objects in the category of totally disconnected Hausdorff spaces over $T$.

Of special interest are so-called {\it corestricted bundles} $\Cal G$ over a profinite space $T$ which are corestricted by a compact subbundle $\Cal U$ of $\Cal G$, i.e.\  $\Cal G$ has the final topology with respect to the inclusions $\Cal U\hookrightarrow\Cal G$ and ${\Cal G}_t\hookrightarrow\Cal G$, $t\in T$, where ${\Cal G}_t$ denotes the fibre of $\Cal G$ at $t$. Alexander Schmidt pointed out that a good way to think of this object is as a ``{\it hedgehog}\,'' having a compact body $\Cal U$ which is surrounded by ``{\it spines}\,'' corresponding to the sets ${\Cal G}_t\backslash {\Cal U}_t$, $t\in T$.

%

In the first section we introduce the notion of (not necessarily compact) bundles and their free products. We investigate bundles endowed with an additional action by a pro-$\c$-group $G$ and their corresponding free products, which are pro-$\c$-$G$ operator groups.
In the second section, we define corestricted bundles and corestricted free products. One main goal of this paper is to study the properties of the functor
\vspace{-.2cm}
$$
F : \hbox{$\cal B${\it undles}}\,\,\lpfeil \,\,\hbox{$\cal P${\it ro}-$\c$-{\it groups}}
$$
which associates to a bundle ${\Cal G}$ over a profinite space $T$ the free pro-$\c$-product ${\freeproduct}_{T} {\cal G}$. This functor commutes with projective limits on the subcategory of compact bundles. We prove that this remains true for certain corestricted $G$-operator bundles. Furthermore, this also holds for projective limits of corestricted bundles over the one-point compactifications of discrete sets.

In the third section, we show that corestricted bundles naturally arise from families of closed subgroups of a pro-$\c$-group. Finally we consider the cohomology of a corestricted free product over the one-point compactification of a discrete set. This yields the number theoretic application we have in mind: Under some conditions there exists a Galois group of a large number field which is a free pro-$p$-product of infinity many  decomposition groups corestricted by inertia groups.

\section{Bundles and free products}

\indent

 Let  $\c$ be a class of finite groups  closed under taking subgroups,
homomorphic images and finite direct products.
A  pro-$\c$-group is a projective limit of groups in $\c$.
In \cite{NSW}, chap.\ IV \S 3, the notion of a bundle of pro-$\c$-groups is introduced. Here we use the notion {\it bundle} in a more general context, but we will follow  partly the presentation of \cite{NSW}; see also \cite{Mel}.

\begin{definition}\label{def-bundle}
Let $T$ be a profinite space, i.e.\ a topological projective limit of finite discrete spaces. A  {\bf bundle of pro-$\c$-groups}\vspace{-.2cm}
$$
{\cal G}=({\cal G},pr,T)
$$

\vspace{-.2cm}\noindent
over $T$ is a group object in the
category of  totally disconnected Hausdorff spaces over $T$
 such that the fibre
over every point of $T$ is a pro-$\c$-group, i.e.\
there are  continuous maps \smallskip {
\Item{} $pr: {\cal G} \rightarrow T$ (the structure map),
\Item{} $m: {\cal G} \times_T {\cal G} \rightarrow {\cal G}$ (the multiplication),
\Item{} $e\,\,: T \rightarrow {\cal G}$ a section to $pr$ (the unit),
\Item{} $\iota\,\,: {\cal G} \rightarrow {\cal G}$  (the inversion),
}

\medskip\noindent
such that the fibre ${\cal G}_t=pr^{-1}(t)$ together with the induced maps
$m_t: {\cal G}_t \times {\cal G}_t \rightarrow {\cal G}_t$, $\iota_t: {\cal G}_t \rightarrow {\cal G}_t$ and
the unit element $e_t=e(t)$ is a pro-$\c$-group for every point $t \in
T$.
\end{definition}

\noindent
Here  ${\cal G} \times_T {\cal G}=\{(x,y)\in {\cal G}\times {\cal G}| pr(x)=pr(y) \}$
denotes the fibre product of ${\cal G}$ and ${\cal G}$ over $T$, which is a
 totally disconnected Hausdorff space and has a structure map $pr\colon {\cal G} \times_T {\cal G}\pfeil T$. Furthermore, the maps $m$ and $\iota$ commute with   the structure maps.

Of special interest are bundles which are compact, i.e.\/ ${\cal G}$ is a profinite space. These are the bundles which were considered in
 \cite{NSW}, and which we now call {\bf compact bundles}.

\medskip
 If $G$ is a pro-$\c$-group and $T$ is a profinite space, then we have the (compact) {\bf constant bundle} $(G\times T,pr, T)$, where $pr$ is the projection $G \times T \rightarrow T$.

\medskip
A {\bf morphism of bundles}
$$
\phi: ({\cal G},pr_{\cal G},T)\rightarrow ({\cal H}, pr_{\cal H},S)
$$
is a pair $\phi_{\cal G} : {\cal G} \rightarrow {\cal H}$,
$\phi_T: T \rightarrow S$ of continuous maps such that the diagram
$$
\xymatrix{
{\Cal G}\ar[r]^{\phi_{\cal G}}\ar[d]_{pr_{\cal G}}&{\Cal H}\ar[d]^{pr_{\cal H}}\\
T\ar[r]_{\phi_T}&S
}
$$
commutes and for every $t \in T$ the associated map
$\phi_t: {\cal G}_t \rightarrow {\cal H}_{\phi_T(t)}$ is a group
homomorphism. We say that $({\cal G},pr_{\cal G},T)$ is a {\bf subbundle} of $({\cal H}, pr_{\cal H},S)$ if $T=S$, $\phi_T=id$ and $\phi_{t}$ is injective for all $t\in T$.
We say that $\phi$ is {\bf fibrewise surjective} if $\phi$ is surjective and $\phi_t$ is surjective for all $t\in T$.

\medskip

A morphism from a bundle $({\cal G},pr,T)$ to a pro-$\c$-group
$G$ is a continuous map $\phi: {\cal G} \rightarrow G$ such that
the induced maps $\phi_t: {\cal G}_t \rightarrow G$ are group
homomorphisms for every $t \in T$.

\subsection{Free pro-$\c$-products}

\indent

Let $({\cal G},pr,T)$ be a bundle of pro-$\c$-groups.

\begin{definition}
The {\bf free pro-$\c$-product} of  $({\cal G},pr,T)$ of
 is  a pro-$\c$-group
$$
G=\freeproduct_{T} {\cal G}
$$
together with a morphism $\omega : {\cal G} \rightarrow G$, which
has the following universal property: for every morphism $f: {\cal
G} \rightarrow H$ from ${\cal G}$ to a pro-$\c$-group $H$ there
exists a unique homomorphism of pro-$\c$-groups $\phi: G
\rightarrow H$ with $f= \phi \circ \omega$.
\end{definition}

As in \cite{NSW}(4.3.6) we have the

\begin{proposition}\label{exsistence}
The free pro-$\c$-product $\freeproduct_{T}{\cal G}$ exists and is
unique up to unique isomorphism.
\end{proposition}
\pagebreak
\demo{Proof:}
Let $D ={\freeproduct}_{t\in T}^{\,d}{\cal G}_t$ be the abstract free product of the family of groups
$({\cal G}_t )_{t \in T}$ and let $\lambda : {\cal G} \rightarrow D$ be
the map which is given on every ${\cal G}_t$ as the natural inclusion of
${\cal G}_t$ into $D$. Then define the free product
$G={\freeproduct}_{T}{\cal G}$ as the completion of $D$ with respect to the
topology which is given by the family of normal subgroups $ N
\subseteq D$ of finite index for which
\Item{(a)} $D/N \in \c$,
\Item{(b)} the composition  of $\lambda$ with the natural
projection $ D \rightarrow D/N$ is continuous.

\smallskip\noindent
It is easily
verified that $G$ has the required universal property. The uniqueness assertion is clear by the universal property.
\enddemo

For the free pro-$\c$-product of the bundle $({\cal G},pr,T)$ we often write
$$
\freeproduct_{t\in T}{\cal G}_t=\freeproduct_T {\cal G}.
$$
Observe that there are different ``free pro-$\c$-products" of a family of pro-$\c$-groups
$\{{\cal G}_t\}_{t\in T}$ depending on the topology on ${\cal G} = \ttcup_{t \in T} {\cal G}_t$.

\medskip
If $\phi: ({\cal G},pr_{\cal G},T)\rightarrow ({\cal H}, pr_{\cal H},S)$ is a morphism of bundles, then by the universal property of the free product the map
$ {\cal G}\pfeil^{\phi_{\cal G}} {\cal H}\pfeil^{\omega}{\freeproduct}_S {\cal H} $ induces a homomorphism\vspace{-.1cm}
$$
\phi_*: \freeproduct_T {\cal G} \rightarrow \freeproduct_S {\cal H},
$$
i.e.\ there is a functor\vspace{-.2cm}
$$
{\Cal B}\lpfeil\hbox{$\cal P${\it ro}-$\c$-{\it groups}}, \hspace{.5cm}{\Cal G}\mapsto \freeproduct_{T} {\cal G},
$$
between the category ${\Cal B}$ of bundles  of pro-$\c$-groups  and the category of  pro-$\c$-groups.

\begin{proposition} \label{zerlegung}
Let $({\cal G},pr,T)$ be a bundle of pro-$\c$-groups and assume that $T=
T_1 \ktcup\cdots\,\ktcup T_k$ is a finite disjoint decomposition of\/ $T$
into non-empty open subsets. Then the following assertions hold.

\Item{(i)}
There is a canonical isomorphism
$$
\freeproduct_{T_1} {\cal G}_1
\freeproductsmall \cdots \freeproductsmall
\freeproduct_{T_k} {\cal G}_k\cong\freeproduct_T {\cal G},
$$
where ${\cal G}_i$ denotes the  bundle $pr^{-1}(T_i)$ over $T_i$
for $i=1,\ldots, k$.

\Item{(ii)}
The canonical homomorphism $\phi_*\colon{\freeproduct}_{T_{i}} {\cal G}_{i}\pfeil{\freeproduct}_T {\cal G} $, which is induced by the bundle morphism
$\phi\colon  {\cal G}_{i}\pfeil {\cal G}$, has a splitting; in particular, $\phi_*$ is injective.
\end{proposition}

\demo{Proof:}
The first assertion follows from the universal property. For the second  consider the map
$\varphi_{{\cal G}}\colon{\cal G}\twoheadrightarrow {\cal G}_i$ which maps ${\cal G}_j$, $j\neq i$, to the unit element of ${\cal G}_{t_0}$ for some fixed $t_0\in T_i$ and which is the identity on ${\cal G}_i$, and the map $\varphi_{T}\colon T\twoheadrightarrow T_i$ which maps $T_j$, $j\neq i$, to  $t_0$ and which is the identity on $T_i$.  Then $(\varphi_{{\cal G}},\varphi_{T})$ is  a morphism of bundles which induces a splitting of $\phi_*$
\enddemo

\begin{corollary} \label{faser-injective1}
Let $({\cal G},pr,T)$ be a bundle, where  $T=T_0\ktcup \{*\}$ is the one-point compactification of the discrete set $T_0$ and ${\cal G}_{*}=\{*\}$. If \ $t_0\in T$, then  the canonical map\vspace{-.2cm}
$$
\omega_ {{\cal G}_{t_0}}: {\cal G}_{t_0}\lpfeil\freeproduct_{t\in T} {\cal G}_t
$$
has a splitting. In particular, it is injective.
\end{corollary}

We do not know whether the map $\omega_ {{\cal G}_{t_0}}$ is injective  for bundles over arbitrary profinite spaces (but see remark 5 in section 3).

\subsection{Compact bundles}

\indent

First we prove a criterion for subsets of a compact bundle for being open.

\begin{lemma}\label{compfamily}
Let ${\Cal G}$ be a compact bundle of pro-$\c$-groups over a profinite space~$T$. For every $t\in T$ let an open subset $W_t$ of ${\Cal G}_t$ be given.  Then the following is equivalent.

\vspace{.2cm}
{
\Item{(i)} For every closed  subset $Y$ of ${\Cal G}$ the set $\{t\in T|\, Y_t\subseteq W_t\}$ is open in $T$.

\vspace{.2cm}
\Item{(ii)} $W=\Tcup_{t\in T}W_t$ is open in  ${\cal G}$.
}
\end{lemma}

\demo{Proof:} Assume that (i) holds. Let
$$
\bar W={\cal G}\backslash W=\Tcup_{t\in T}\bar W_t,\hspace{.5cm} \bar W_t={\Cal G}_t\backslash W_t.
$$
Let $g\in  W_{t_0}\subseteq W$. Since $\bar W_{t_0}$ is closed in ${\Cal G}_{t_0}$ (and so in ${\Cal G}$), there are open subsets  $R,Q$ of ${\Cal G}$ such that $\bar W_{t_0}\subseteq Q$, $g\in R$ and $R\cap Q= \varnothing$. Since $\bar Q={\cal G}\backslash Q$ is closed in ${\Cal G}$, the set
$$
S:=\{t\in T| \bar W_t\subseteq Q\}=\{t\in T| \bar Q_t \subseteq W_t\}
$$
is open in $T$ by assumption and $t_0\in S$. Thus $V=pr_{\Cal U}^{-1}(S)\cap R$  is an open neighborhood of $g$. Furthermore, if $t\in S$, then
 $V_t\subseteq R_t\subseteq \bar Q_t\subseteq W_t$.
Thus $V$ is contained in $W$. It follows that  $W$ is open.

Now assume that (ii) holds. Then\vspace{-.1cm}
$$
\{t\in T|\, Y_t\subseteq W_t\}=T\backslash pr(Y\cap \bar W).
$$
Since $(Y\cap \bar W)$ is closed in ${\cal G}$, hence compact, $pr(Y\cap \bar W)$ is a compact subset of $T$, and so $T\backslash pr(Y\cap \bar W)$ is open.
\enddemo

Let $G$ be a pro-$\c$-group and let $(G_t)_{t \in T}$
be a family of closed subgroups of $G$ indexed by the points of a profinite space
$T$.  We say  that
$(G_t)_{t \in T}$ is a {\bf continuous family}, if and only if
for every open neighborhood  $V$ of the identity of  $G$ the set
$T(V):=\{t \in T\,|\, G_t \subseteq V \}$ is open in $T$.

\pagebreak

The following proposition shows that the concepts  of compact bundles and continuous families coincide.

\begin{proposition}\label{family}

\noindent\smallskip
\Item{(i)}
The set
$$
{\cal G} =\{(g,t)\in G \times T\,|\, g \in G_t \}
$$
equipped with the induced topology of the constant bundle $(G\times T,pr,T) $ is a compact bundle over $T$.

\noindent\smallskip
\Item{(ii)} If ${\cal H}$ is a compact bundle, then ${\cal H}$ is the bundle associated to the family $({\cal H}_t)_{t \in T}$ of pro-$\c$-groups given by the fibres of ${\cal H}$, where each fibre is considered as closed subgroup of ${\freeproduct}_T{\cal H}$.

\end{proposition}

\demo{Proof:} For the  first assertion see \cite{NSW} (4.3.3). In order to prove the second, we   first remark that  a fibre  ${\cal H}_t$ is a closed subgroup of the free pro-$\c$-product ${\freeproduct}_T{\cal H}$, see \cite{NSW} (4.3.12)(i), thus\vspace{-.2cm}
$$
{\cal H}=\{(h,t)\in  \freeproduct_T{\cal H}\times T\,|\, h \in {\cal H}_t\}.
$$
 Now it is easy to see that the family $({\cal H}_t)_{t \in T}$ is continuous: Let $W$ be an open neighborhood of the unit of ${\freeproduct}_T{\cal H}$. Then, using lemma (\ref{compfamily}), the set $T(W)=\{t \in T\,|\, {\cal H}_t \subseteq W \}=\{t \in T\,|\, {\cal H}_t \subseteq (W\cap {\Cal H})_t \}$ is open in $T$.
\enddemo

\subsection{Bundles with operators}

\indent
\vspace{-.2cm}
\begin{definition} Let $G$ be a pro-$\c$-group. A  bundle  $({\cal G},pr, T)$ of  pro-$\c$-groups is called {\bf $G$-bundle} if  $G$ acts continuously on ${\cal G}$, i.e.\ there is a commutative diagram
$$
\xymatrix{
G\ar@<-3.2ex>@{=}[d]\,\,\times \,\, {\Cal G}\ar[r]\ar@{->>}@<3.ex>[d]^-{pr}&{\Cal G}\ar@{->>}[d]^-{pr}& (\sigma,x)\mapsto \sigma x,\\
G\,\,\times \,\,T\ar[r]&T& (\sigma,t)\mapsto \sigma t,
}
$$
of continuous maps such that

\smallskip
\Item{(i)} for all $\sigma\in G$ the map ${\cal G}\pfeil {\cal G},\ x\mapsto \sigma x$, is a morphism of bundles,

\smallskip
\Item{(ii)} $(\sigma\tau)(x)=\sigma(\tau(x))$ and $ex=x$ for all $\sigma,\tau\in G$ and all $x\in {\cal G}$.
\end{definition}

\begin{definition}
Let  $({\cal G},pr, T)$ and $({\cal H},pr, S)$ be $G$-bundles of  pro-$\c$-groups and let $\phi\colon {\cal G}\pfeil {\cal H} $ be a morphism of bundles. Then $\phi$ is called {\bf $G$-invariant}, if the diagrams\vspace{-.2cm}
$$
\xymatrix{
G\ar@<-3.2ex>@{=}[d]\,\,\times \,\, {\Cal G}\ar[r]\ar@<3.ex>[d]^-{\phi_{\Cal G}}&{\Cal G}\ar[d]^-{\phi_{\cal G}}&&
G\ar@<-3.2ex>@{=}[d] \,\,\times \,\,T\ar[r]\ar@<3.ex>[d]^-{\phi_T}&T\ar[d]^-{\phi_T}\\
G\,\,\times \,\, {\Cal H}\ar[r]&{\Cal H}&&G \,\,\times \,\,S\ar[r]&S
}
$$
commute. The morphism $\phi$ is called {\bf $G$-transitive}, if $\phi$ is $G$-invariant and  for
$t,t'\in\phi_T^{-1}(s)$, $s\in S$, there exists $\sigma\in G$  such that $\sigma t=t'$.
\end{definition}

Let $({\cal G},pr, T)$ be a $G$-bundle and let $\sigma\in G$. Then
$$
\sigma {\cal G}_t={\cal G}_{\sigma t}\hspace{.5cm}\hbox{for $t\in T$}.
$$

\noindent
Furthermore, the homeomorphism
$
\sigma \colon {\cal G}\liso {\cal G}
$
induces an  automorphism
$$
\sigma : \freeproduct_{T} {\cal G}\liso \freeproduct_{T} {\cal G}
$$
of pro-$\c$-groups, i.e.\ ${\freeproduct}_{T} {\cal G}$ is a pro-$\c$-$G$ operator group, see \cite{NSW} (4.3.8).  Thus ${\freeproduct}_{T} {\cal G}$ possesses a system of neighborhoods of the identity consisting of open  $G$-invariant normal subgroups, see \cite{P} theorem 17:

If $U$ is an open subgroup of ${\freeproduct}_{T} {\cal G}$, then
$
\tilde U :=\bigcap_{\sigma\in G}\sigma U
$
is a $G$-invariant  subgroup of $ {\freeproduct}_{T} {\cal G}$. We will show that $\tilde U$ is open. Since $\sigma e=e\in U$, there exist open neighborhoods $V_\sigma$ and $W_\sigma$ of $e\in U$ and $\sigma\in G$, respectively, such that $W_\sigma V_\sigma\subseteq U$. Since $\bigcup_{\sigma}W_\sigma=G$ and $G$ is compact, we find a finite covering $\{W_{\sigma_1},\dots,W_{\sigma_n}\}$ of $G$. Let $V= \bigcap_{i=1}^n V_{\sigma_i}$, then $\sigma V\subseteq \tilde U$ for every $\sigma\in G$ and consequently for every $x\in \tilde U$ we have $Vx\subseteq \tilde U$.

\medskip

Since the map $\lambda\colon {\cal G}\pfeil {\freeproduct}_{T}^d {\cal G}$ is $G$-invariant, where the action of $G$ on ${\freeproduct}_{T}^d  {\cal G}$
is defined in the obvious way, it follows that
$$
\freeproduct_{T} {\cal G}=\varprojlim_N({\freeproduct_T}^d {\cal G})/N,
$$

\vspace{-.2cm}
where $N$ runs through all  $G$-invariant normal subgroups of ${\freeproduct}_{T}^d {\cal G}$ such that \linebreak $({\freeproduct}_{T}^d {\cal G})/N\in \c$ and the map
${\cal G}\pfeil {\freeproduct}_{T}^d {\cal G}\twoheadrightarrow
 ({\freeproduct}_{T}^d {\cal G})/N$  is continuous.

\subsection{Projective limits of bundles}

Let $I$ be a directed set and let $\{{\cal G}_i,T_i,\phi_{ij}\}_I$ be a projective system of bundles of pro-$\c$-groups with transition morphisms
 $$
 \phi_{ij}: {\cal G}_i\to {\cal G}_j,\ i\ge j.
 $$
 A  bundle $({\cal G},pr,T)$ of pro-$\c$-groups together with morphisms of  bundles $\phi_i\colon {\cal G} \to {\cal G}_i$ compatible with $\phi_{ij}$ is called {\bf projective limit} of $\{{\cal G}_i,T_i,\phi_{ij}\}_I$, if it satisfies the usual universal property. The projective limit  exists and is unique up to isomorphism and we write
$$
({\cal G},pr,T) = \varprojlim_{i \in I} ({\cal G}_i,pr_i,T_i).
$$
Indeed, the uniqueness follows from the universal property and for the existence we consider the Cartesian product $\tilde{{\cal G}}:=\prod_{i\in I} {\cal G}_i$ which, being endowed with the product topology, is a bundle over $T=\varprojlim T_i$. Then the subbundle
\[
 {\cal G}:= \{(g_i)\in \tilde{{\cal G}}\ |\ \phi_{ij} (g_i) = g_j\}
\]
 satisfies the universal property of the projective limit and its fibres are given by
\[
 {\cal G}_t =  \varprojlim_{i \in I} {\cal G}_{i, t_i}
\]
for $t=(t_i)_{i\in I}\in \varprojlim T_i$. Clearly, the transition morphisms $\phi_{ij}$ give rise to homomorphisms $\freeproduct_{T_i} {\cal G}_i\pfeil \freeproduct_{T_j} {\cal G}_j,\ i\ge j$. It is a natural question whether the induced homomorphism\vspace{-.2cm}
\[
 \freeproduct_T {\cal G}\lpfeil \varprojlim_{i \in I} \freeproduct_{T_i} {\cal G}_i
\]
is an isomorphism, i.e.\ whether free products commute with projective limits. This holds if the bundles ${\cal G}_i$ and hence ${\cal G}$ are compact, see \cite{NSW}(4.3.6). Unfortunately, the argument cannot be carried over to arbitrary bundles. However, under additional conditions for the maps ${\cal G}\to {\cal G}_i$, we can prove the following criterion:

\begin{proposition}\label{projlim1} Let $G$ be a pro-$\c$-group and
let $\{{\cal G}_i,T_i,\phi_{ij}\}_I$ be a projective system of  $G$-bundles of pro-$\c$-groups with $G$-invariant transition maps $\phi_{ij}$. Then the  bundle
$$
{\cal G} = \varprojlim_{i \in I} {\cal G}_i
$$
is a $G$-bundle over  $T={\varprojlim} T_i$. Assume further that the following holds:

{
\smallskip
\Item{(i)} For all $i\in I$, the canonical map $\phi_{i}\colon {\cal G}\pfeil {\cal G}_i$ is surjective.

\smallskip
\Item{(ii)} Let $W$ be a closed resp.\ open  subset of ${\cal G}$ which is invariant under  an  open  subgroup $M$ of $G$. Then $\phi_i(W)$ is closed resp.\ open  in
$ {\cal G}_i$ for all $i\in I$.

}

\noindent\medskip
Then
$$
\freeproduct_T {\cal G} =
 \varprojlim_{i \in I} \freeproduct_{T_i} {\cal G}_i.
$$
\end{proposition}

\demo{Proof:}
The fact that ${\cal G}$ is a $G$-bundle over $T$ is obvious. Let $D ={\freeproduct}_{t\in T}^{\,d}{\cal G}_t$ be the abstract free product of the family of groups
${\cal G}_t$
 and
$D_i ={\freeproduct}_{t_i\in T_i}^{\,d}{\cal G}_{i,t_i}$ for $i\in I$.
Using (i), we see that the canonical maps
$$
\varphi_i\colon D\surjr{}D_i,\hspace{.5cm} i\in I,
$$
are surjective. Furthermore the action of $G$ on ${\cal G}$ induces a $G$-action on $D_i$ and $D$, and $\varphi_i$ is $G$-invariant. For a fixed $i\in I$, the correspondence $N\mapsto N_i:=\varphi_i(N)$ induces a bijection between the set of $G$-invariant normal subgroups $N$ of $D$ such that $N\supseteq \ker \varphi_i$ and $\ D/N\in\c$, and the set of normal subgroups $N_i$ of $D_i$ such that $D_i/N_i\in\c$. For all such $N$ and $N_i$, we have the isomorphism
$$
D/N\liso D_i/N_i
$$
and the commutative diagram\vspace{-.2cm}
$$ \xyalign
\xymatrix{
{\Cal G}_i\ar[r]^-{\lambda_i}& D_i \ar@{->>}[r]^-{f_i}&D_i/N_i\\
{\Cal G}\ar[r]^-{\lambda}\ar@{->>}[u]^-{\phi_i}&D\ar@{->>}[u]^-{\varphi_i}\ar@{->>}[r]^-{f}&D/N\,.\ar[u]^-\sim\\
}
$$

\noindent
Obviously, if $f_i\lambda_i$ is continuous, then $f\lambda$ is continuous. Conversely, assume that $f\lambda$ is continuous and let $a\in D_i/N_i$. Then
$W=(f\lambda)^{-1}(a)$ is open and closed in ${\Cal G}$ and invariant under any open subgroup of $G$ acting trivially on the finite group $D/N$. Furthermore, since $\phi_i$ is surjective, $\phi_i(W)=(f_i\lambda_i)^{-1}(a)$.
Using (ii),  $\phi_i(W)$ is closed resp.\ open in ${\Cal G}_i$ . This shows that $f_i\lambda_i$ is continuous.
\smallskip

Now by the general theory of pro-$\c$-$G$ operator groups presented in the previous section, we have $\freeproduct_T {\cal G} = \varprojlim_{N}D/N$ where $N$ runs through the $G$-invariant normal subgroups of $D$ such that $D/N\in \c$ and $f\lambda$ is continuous. Analogously, $\freeproduct_T {\cal G}_i = \varprojlim_{N_i}D_i/N_i$ where $N_i$ runs through the $G$-invariant normal subgroups of $D_i$ such that $D_i/N_i\in \c$ and $f_i\lambda_i$ is continuous. Noting that any normal subgroup $N$ of $D$ of finite index contains $\ker \varphi_i$ for some $i\in I$, the above considerations imply
\[
 \freeproduct_T {\cal G} = \varprojlim_{N}D/N = \varprojlim_{i\in I}\big( \varprojlim_{N\supseteq \ker(\varphi_i)} D/N\big) = \varprojlim_{i\in I}\big( \varprojlim_{N_i} D_i/N_i\big) = \varprojlim_{i\in I} \freeproduct_{T_i} {\cal G}_i.
\]
\enddemo

\noindent
{\bf Remark:} In section 2 we will work with projective limits which are formed in the subcategory of {\bf corestricted bundles}. As a set, the projective limit ${\cal G}=\varprojlim {\cal G}_i$ is given as above by ${\cal G}:= \{(g_i)\in \prod {\cal G}_i\ |\ \phi_{ij} (g_i) = g_j\}$, however it has to be endowed with a topology which is finer than the product topology. It is important to remark that (\ref{projlim1}) remains valid as long as the projections ${\cal G}\pfeil {\cal G}_i$ are continuous and condition (ii) holds with respect to this finer topology.

\section{Corestricted bundles and their free products}

\indent

Now we introduce the notion of a corestricted bundle of pro-$\c$-groups with
respect to a compact bundle ${\cal U}$  of pro-$\c$-groups.

\begin{definition}\label{def-cor-bundle} Let $({\cal G},pr_{{\cal G}},T)$ be a bundle of pro-$\c$-groups. We say that ${\cal G}$ is
 {\bf corestricted } with respect to a compact bundle $({\cal U},pr_{{\cal U}},T)$ over $T$ if there exists an injective  morphism of bundles
$\phi\colon {\cal U}\hookrightarrow {\cal G}$ over $T$, i.e.\
the diagram\vspace{-.2cm}
$$
\xymatrix{
\cal U\ar@{^(->}[rr]^{\phi}\ar[dr]_{pr_{\cal U}}&& \Cal G\ar[dl]^{pr_{{\cal G}}}\\
&T&\\
}
$$
commutes, and
${\cal G}$ is equipped with the final topology with respect to the family of inclusions
$\{\phi\colon {\cal U}\hookrightarrow {\cal G},\,\, {\cal G}_t\hookrightarrow {\cal G}, t\in T \}$.
We denote the bundle
${\cal G}$ by $({\cal G},{\cal U})$ if it is corestricted with respect to ${\cal U}$.

\medskip
The {\bf corestricted free  pro-$\c$-product} of the corestricted bundle $({\cal G},{\cal U})$
is the free pro-$\c$-product ${\freeproduct}_{T} ({\cal G},{\cal U})$,
and we write
$$
\freeproduct_{t\in T} ({\cal G}_t,{\cal U}_t)=\freeproduct_{T} ({\cal G},{\cal U}).
$$
\end{definition}

\vspace{.5cm}\noindent
{\bf Remarks:}
%
1. The topology on ${\cal G}_t$ induced by the topology of $({\cal G},{\cal U})$ is the pro-$\c$ topology of ${\cal G}_t$: if $V_0$ is an open subgroup of the pro-$\c$-group ${\cal G}_{t_0}$, then ${\cal G}_{t_0}\backslash V_0$ is closed in ${\cal G}_{t_0}$ and so in ${\cal G}$ since every fiber of ${\cal G}$ is closed. It follows that $V=\ttcup_{t \neq t_0} {\cal G}_t\ktcup V_0$ is open in ${\cal G}$, and $V\cap {\cal G}_{t_0}=V_0$.

\noindent
2. If $T$ is the one-point compactification of a discrete set $T_0$, then
the definition of
the corestricted free pro-$\c$-product of the corestricted bundle $({\cal G},{\cal U})$ over
the profinite space $T$ coincides
with that given in \cite{N}\,\S 2.  If $T$ is finite, then $({\cal G},{\cal U})$ is a  compact bundle independent of ${\cal U}$.

\noindent
3. We have two extreme  cases.
If ${\cal U}_t={\cal G}_t$ for all  $t\in T$, then $({\cal G},{\cal U})$ is the
equal to the compact bundle ${\cal U}$. If ${\cal U}_t=\{1\}$ for all  $t\in T$, i.e.\
${\cal U}$ is the trivial bundle $\it 1$, then the topology of  $({\cal G},{\it 1})$ is given by the  sets
$V\subseteq {\cal G}$ such that
$V_t=V\cap {\cal G}_t$ is an open subgroup of ${\cal G}_t$ for all $t\in T$ and  $S=\{t\in T|1\in V_t\}$ is open in $T$.
The corestricted free  pro-$\c$-product
$$
\ffreeproduct_{T} ({\cal G},{\it 1})
$$
of the
corestricted bundle  $({\cal G},{\it 1})$ is  sometimes called the {\it unrestricted
free  pro-$\c$-product} of ${\cal G}$.

\bigskip
\begin{definition}\label{def-cor-morphism}
Let $T,S$ be profinite spaces and $({\cal G}, {\cal U}),\ ({\cal H}, {\cal V})$ be corestricted bundles of pro-$\c$-groups over $T$ and $S$, respectively. A morphism of bundles
\[
 \phi: ({\cal G}, {\cal U}) \lpfeil ({\cal H}, {\cal V})
\]
is called {\bf morphism of corestricted bundles} if $\phi({\cal U})\subseteq {\cal V}$. Furthermore, we say that $\phi$ is {\bf strict} if $\phi^{-1}({\cal V})={\cal U}$.
\end{definition}

As a special case of morphisms of corestricted bundle, we consider morphisms $(\cal G, {\cal U}\,')\to (\cal G,\cal U)$ where ${\cal U}\,'\subseteq {\cal U}$ are compact subbundles. In this situation, we have the following

\begin{proposition} \label{subbundles}
 Let $(\cal G, \cal U)$ and $(\cal G,{\cal U}\,')$ be corestricted bundles of pro-$\c$-groups with respect to compact subbundles $\cal U$ and ${\cal U}\,'$, respectively, with continuous  inclusion ${\cal U}\,'\subseteq \cal U$. Then the morphism of corestricted bundles $id: (\cal G, {\cal U}\,')\to (\cal G,\cal U)$ induces a canonical surjection
$$
 \freeproduct_{T} ({\cal G}, {{\cal U}\,'}) \surjr{} \freeproduct_{T} (\cal G, \cal U)
$$
and there is an isomorphism
$$
   \varprojlim_{N} (\freeproduct_T ({\cal G},{{\cal U}\,'}))/{N} \liso{} \freeproduct_T ({\cal G}, {\cal U}).
$$
Here $N$ runs through the open normal subgroups of ${\freeproduct}_T ({\cal G},{{\cal U}\,'})$ such that for every $a\in {\freeproduct}_T ({\cal G},{{\cal U}\,'})$ the preimage of $aN$ under the map
$$
 ({\cal G},{\cal U})  \stackrel{id}{\hbox{$\lpfeil$}} ({\cal G},{{\cal U}\,'}) \stackrel{\omega}{\hbox{$\lpfeil$}} \freeproduct_{T} (\cal G,{\cal U}\,')
$$
is open in $(\cal G,\cal U)$, i.e.
 \Item{(i)} $\omega^{-1} (aN)\cap \cal U$ is open in $\cal U$,
 \Item{(ii)} $\omega^{-1} (aN) \cap {\cal G}_t$ is open in ${\cal G}_t$ for all $t\in T$.
\end{proposition}
\demo{Proof:}
This follows directly from the universal property of the free product. In fact, ${\freeproduct}_T ({\cal G}, {\cal U})$ satisfies the property of the completion of ${\freeproduct}_T ({\cal G},{{\cal U}\,'})$ with respect to the family of open normal subgroups $N$ satisfying (i) and (ii). Hence the homomorphism ${\freeproduct}_{T} ({\cal G}, {{\cal U}\,'}) \surjr{} {\freeproduct}_{T} (\cal G, \cal U)$ induces the isomorphism $\varprojlim_{N} ({\freeproduct}_T ({\cal G},{{\cal U}\,'}))/{N} \liso{} {\freeproduct}_T ({\cal G}, {\cal U})$.
\enddemo

In particular, it follows that any corestricted free product can be recovered from the unrestricted free product ${\freeproduct}_T ({\cal G},{\it 1})$.

\subsection{Quotients of corestricted bundles}

\indent

Let $(\cal G, \cal U)$ be  a corestricted bundle of pro-$\c$-groups and let ${\cal V}$ be a closed subbundle of ${\cal G}$ such that ${\cal V}_t$ is a normal subgroup of ${\cal G}_t$ for all $t\in T$. Set $\overline{\cal V}:={\cal V}\cap {\cal U}$ and let $\tilde{\cal U}:= \ttcup_{t\in T} {\cal U}_t/\overline{\cal V}_t$ be endowed with the quotient topology with respect to the canonical surjection ${\cal U}\twoheadrightarrow\tilde{\cal U}$. Then $\tilde{\cal U}$ is a compact bundle over $T$ and we define the {\bf quotient bundle} $(\tilde{\cal G}, \tilde{\cal U}):=({\cal G}/{\cal V},{\cal U}/\overline{\cal V})$ of $({\cal G},{\cal U})$ with respect to ${\cal V}$ as the set
$$
 (\tilde{\cal G}, \tilde{\cal U})=\Tcup_{t\in T} {\cal G}_t/{\cal V}_t
$$
endowed with the final topology with respect to the inclusions  ${\cal G}_t/{\cal V}_t\hookrightarrow (\tilde{\cal G}, \tilde{\cal U})$, $t\in T$, and
$\tilde{\cal U}\hookrightarrow (\tilde{\cal G}, \tilde{\cal U})$. Thus $(\tilde{\cal G}, \tilde{\cal U})$ is a corestricted bundle of pro-$\c$-groups and there is a canonical surjective morphism of corestricted bundles\vspace{-.2cm}
$$
\phi : ({\cal G}, {\cal U}) \surjr{} (\tilde{\cal G},\tilde{\cal U})
$$
which induces a surjection\vspace{-.3cm}
\[
 \freeproduct_T ({\cal G},{\cal U}) \surjr{} \freeproduct_T (\tilde{\cal G},\tilde{\cal U}).
\]

\noindent
Indeed, $\phi$ is continuous: Let $\tilde V$ be open in $(\tilde{\cal G}, \tilde{\cal U})$, and so $\tilde V \cap \tilde {\Cal U}$ is open in $\tilde {\Cal U}$.  Then the equality
$
\phi^{-1}(\tilde V)\cap  {\Cal U}=\phi^{-1}(\tilde V \cap \tilde {\Cal U})\cap {\Cal U}=\phi_{|{\Cal U}}^{-1}(\tilde V \cap \tilde {\Cal U})
$
shows that $ \phi^{-1}(\tilde V)\cap  {\Cal U}$ is open in ${\Cal U}$. With an analogous argument for the fibres it follows that  $\phi^{-1}(\tilde V)$ is open in $({\cal G},{\cal U})$.

\bigskip\noindent
{\bf Remark:}\label{quotientbundle}
Let $(\cal G, \cal U)$ be a corestricted bundle of pro-$\c$-groups such that ${\cal U}_t$ is a normal subgroup of ${\cal G}_t$ for all $t\in T$.
Then the topology of the quotient bundle $({\cal G}/{\cal U},{\it 1})$, where ${\it 1}$ is the compact unit bundle $\ttcup_{t\in T} \{{\it 1}_t\}$, can be described as follows:
A subset $V\subseteq ({\cal G}/{\cal U},{\it 1} )$ is open if and only if

{
\smallskip
\Item{(i)} $V\cap {\cal G}_t/{\cal U}_t$ is open in ${\cal G}_t/{\cal U}_t$ for all $t\in T$,
\Item{(ii)}  $T(V)=\{t\in T\ |\ 1_t\in V\}$ is open in $T$ where $1_t$ denotes the unit in ${\cal G}_t/{\cal U}_t$.
}

\noindent
In fact,  a subset $V$ of  $({\cal G}/{\cal U},{\it 1})$ is open if and only if $\pi^{-1}(V)$ is open in $({\cal G},{\cal U})$, i.e.\ if and only if $\pi^{-1}(V)\cap {\cal U}$ is open in $\cal U$ and  $\pi^{-1}(V)\cap {\cal G}_t$ is open in ${\cal G}_t$ for any $t\in T$. The last statement is equivalent to (i), and since
\[
 \pi^{-1}(V)\cap {\cal U} = \Tcup_{t\in T(V)} {\cal U}_{\,t}
\]
the first statement is equivalent to the assertion that
$T(V)$ is open in $T$.

\indent

Of particular interest are quotient bundles which occur when passing to $\tilde\c$-completions, where $\tilde\c$ is a class of finite groups which is closed under taking subgroups,
homomorphic images and finite direct products and contained in the class $\c$. If $G$ is a pro-$\c$-group, then $G(\tilde\c)$ denotes its maximal pro-$\tilde\c$ factor group.

\begin{proposition}\label{changing} Let $({\cal G}, {\cal U})$ be a corestricted bundle of pro-$\c$-groups over $T$.
Assume that either

{
\Item{} ${\cal G}={\cal U}$ is a compact bundle or
\Item{} $T=T_0\cup\{*\}$ is the one-point compactification of a discrete set $T_0$

}

\noindent
{\rm (}and $\Cal G_*= \Cal U_*=\{*\}${\rm )}. Then, with the notion as above, we have a surjective morphism of  corestricted bundles\pagebreak
$$
\phi :  ({\cal G}, {\cal U})\twoheadrightarrow (\tilde{\cal G}, \tilde{\cal U}).
$$
This morphism induces  an isomorphism
$$
\big( \freeproduct_T ({\cal G},{\cal U})\big)(\tilde\c)\cong\tilde{\freeproduct_T} (\tilde{\cal G},\tilde{\cal U}),
$$
where $\tilde\freeproduct$ denotes the free product  the category of pro-$\tilde\c$-groups.
\end{proposition}

\demo{Proof:} Denote by ${\cal V}_t$ the kernel of the homomorphism $\rho_t\colon {\cal G}_t \pfeil \tilde{\cal G}_t,\ t\in T$. Once we have shown that  ${\cal V}:=\ttcup_{t\in T} {\cal V}_t$ is a closed subbundle of $({\cal G},{\cal U})$, we  proceed as follows:

\smallskip\noindent
The bundle ${\Cal U}\cap {\cal V}$ is closed in
$({\cal G},{\cal U})$ and  therefore $(\tilde{\cal G}, \tilde{\cal U})=({\cal G}/{\cal V},{\cal U}/({\cal V}\cap {\cal U}))$ is the quotient bundle of $({\cal G},{\cal U})$ with respect to ${\cal V}$ and the canonical surjection $\phi \colon ({\cal G}, {\cal U})\twoheadrightarrow (\tilde{\cal G}, \tilde{\cal U})$ is a morphism of corestricted bundles. Furthermore, $\phi$ induces a surjection
$$
\phi_*:\big( \freeproduct_T ({\cal G},{\cal U})\big)(\tilde\c)
\twoheadrightarrow\tilde{\freeproduct_T} (\tilde{\cal G},\tilde{\cal U}).
$$
On the other hand, the morphism
$\varphi\colon({\cal G},{\cal U})\pfeil^{\omega}{\freeproduct}_T ({\cal G},{\cal U})\twoheadrightarrow \big( \freeproduct_T ({\cal G},{\cal U})\big)(\tilde\c)$
induces a map
$$
\tilde\varphi\colon(\tilde{\cal G},\tilde{\cal U})\pfeil\big( {\freeproduct}_T ({\cal G},{\cal U})\big)(\tilde\c).
$$
We will show that $\tilde\varphi$ is continuous. Let $W$ be open in $ \big( {\freeproduct}_T ({\cal G},{\cal U})\big)(\tilde\c)$ and let
$\tilde V=\tilde\varphi^{-1}(W)$. Then $V:=\varphi^{-1}(W)$ is open in $({\cal G},{\cal U})$, and so $V\cap  {\Cal U}$ is open in ${\Cal U}$. Since $\tilde{\cal U}$ has the quotient topology with respect to the surjection ${\cal U}\twoheadrightarrow\tilde{\cal U}$, the equality
$$
V\cap  {\Cal U}=\phi^{-1}(\tilde V \cap \tilde {\Cal U})\cap {\Cal U}=\phi_{|{\Cal U}}^{-1}(\tilde V \cap \tilde {\Cal U}),
$$
shows that $\tilde V \cap \tilde {\Cal U}$ 
is open in $\tilde{\Cal U}$. Using an analogous argument for the intersection of $\tilde V$ with a fibre, it follows that $\tilde V$ is open in
$(\tilde{\cal G},\tilde{\cal U})$.

Thus we get a homomorphism $\tilde{\freeproduct}_T (\tilde{\cal G},\tilde{\cal U})\pfeil  \big( {\freeproduct}_T ({\cal G},{\cal U})\big)(\tilde\c)$ which is inverse to $\phi_*$. Now we will show that under our assumptions ${\cal V}$ is a closed subbundle.

If ${\cal G}={\cal U}$ is a compact bundle, then ${\cal G}$ is projective limit of finite bundles ${\cal G}_\lambda$, see \cite{NSW} (4.3.10). Since
$\tilde{\cal G}_\lambda$ is obviously a (compact) bundle of pro-$\tilde \c$-groups, it follows that $\tilde{\cal G}=\varprojlim_\lambda\tilde{\cal G}_\lambda$ is a compact bundle of pro-$\tilde \c$-groups and  ${\cal V}$ is a closed subbundle of ${\cal G}$.

If $T=T_0\cup\{*\}$ is the one-point compactification of a discrete set $T_0$, then
the fibres $\Cal U_t$, $t\in T_0$, are open in $\Cal U$. Thus the set
$\ttcup_{t\in T_0}\Cal U_t\backslash \Cal V_t$ is open in $\Cal U$, hence ${\cal V}=\ttcup_{t\in T}\Cal V_t$ is closed.
\enddemo

If $\a$ is the class of finite abelian groups, then  the free pro-$\a$-product and the free pro-$(\c\cap\a)$-product  of abelian pro-$\c$-groups coincide. Thus we obtain the following

\begin{corollary}\label{abelian} Let $({\cal G},{\cal U})$ be  a corestricted bundle of pro-$\c$-groups over $T$,
where is the one-point compactification of a discrete set $T_0$ and $\Cal G_*= \Cal U_*=\{*\}$. Then there is  a canonical isomorphism
$$
\big( \freeproduct_{t\in T} ({\cal G}_t,{\cal U}_t)\big)^{ab}\cong\tilde{\freeproduct_{t\in T}} ({\cal G}_t^{ab},\tilde{\cal U}_t),
$$
where $\tilde{\cal U}_t={\cal U}_t[{\cal G}_t,{\cal G}_t]/[{\cal G}_t,{\cal G}_t]$ and $\tilde\freeproduct$ denotes the free product in the category of pro-abelian-groups.
\end{corollary}

\subsection{Projective limits of corestricted bundles}
\label{coreslimitsection}

\indent

If  $\{({\cal G}_i, {\cal U}_i),T_i,\phi_{ij}\}_I$ is a projective system of corestricted bundles of pro-$\c$-groups, i.e.\ the transition morphisms
 $\phi_{ij}$  are morphisms of corestricted bundles, then the projective limit of this system as defined in section $1$ need not to be a corestricted bundle. The reason is that the topology on the limit is too coarse. Nevertheless, the category of corestricted bundles is closed under projective limits: We endow the Cartesian product $\tilde{{\cal G}}=\prod_{i\in I} {\cal G}_i$ with the final topology with respect to the inclusions
$
 \prod_{i\in I} {\cal G}_{i,t_i} \hookrightarrow \tilde{{\cal G}}
$
 for any $(t_i)\in \prod T_i$  and $\tilde{{\cal U}} \hookrightarrow \tilde{{\cal G}}$ where $\tilde{{\cal U}}:=\prod_{i\in I} {\cal U}_i$ is endowed with the product topology. Then  $\tilde{{\cal G}}$ is a corestricted bundle with respect to the compact subbundle $\tilde{{\cal U}}$. Furthermore,\vspace{-.1cm}
\[
 {\cal G}:= \{(g_i)\in \tilde{{\cal G}}\ |\ \phi_{ij} (g_i) = g_j\}
\]
is a corestricted bundle with respect to the compact subbundle
\[
  {\cal U}:= \{(u_i)\in \tilde{{\cal G}}\ |\ \phi_{ij} (u_i) = u_j\}
\]
satisfies the universal property of the projective limit. We  write
\[
 ({\cal G}, {\cal U})= \itvarprojlim{i \in I}({\cal G}_i, {\cal U}_i).
\]

\begin{theorem}\label{corestricted-proj-lim} Let $G$ be a pro-$\c$-group and
let $\{({\cal G}_i, {\cal U}_i),T_i,\phi_{ij}\}_I$ be a projective system of corestricted $G$-bundles of pro-$\c$-groups with $G$-invariant transition maps $\phi_{ij}$. Then the corestricted bundle
$$
({\cal G},{\cal U}) = \itvarprojlim{i \in I} ({\cal G}_i, {\cal U}_i)
$$
over $T={\varprojlim} T_i$ is a $G$-bundle. Assume that

{
\smallskip
\Item{(i)} the  action of $G$  on each $T_i$  factors through some  finite quotient of $G$,

\smallskip
\Item{(ii)} the transition maps $\phi_{ij}$ are  surjective strict morphisms and

\smallskip
\Item{(iii)}the  maps $\phi_{i}\colon ({\cal G},{\cal U})\twoheadrightarrow ({\cal G}_i, {\cal U}_i)$ are $G$-transitive.

}

\smallskip\noindent
Then
$$
\freeproduct_T ({\cal G}, {\cal U}) =
 \varprojlim_{i \in I} \freeproduct_{T_i} ({\cal G}_i, {\cal U}_i).
$$

\end{theorem}

 Using the following lemma, the above result follows from (\ref{projlim1}). Note that (\ref{projlim1}) holds with respect to the topology on the corestricted projective limit bundle $({\cal G}, {\cal U})$, cf.\ the remark at the end of section 1.

\begin{lemma}\label{closed-mapping-lemma} Let $G$ be a pro-$\c$-group acting on the corestricted bundles
$({\cal G}, {\cal U})$ and $({\cal H}, {\cal V})$  of pro-$\c$ groups over profinite spaces $T$ and $S$, respectively. Assume that the action of $G$ on $S$ factors through a finite factor group $G/N$. Let
$$
\phi\colon ({\cal G}, {\cal U}) \surjr{} ({\cal H}, {\cal V})
$$
be a $G$-transitive,
surjective, strict morphism. Then the following holds.

\smallskip
\Item{(i)} Every fibre ${\cal H}_s$, $s\in S$, is a finite union of images $\phi({\cal G}_{t})$ of fibres of ${\cal G}$.

\smallskip
\Item{(ii)}
Let $W$ be a closed subset of $({\cal G}, {\cal U})$ which is invariant under an open  subgroup $M$ of $G$. Then $\phi(W)$ is closed in $({\cal H}, {\cal V})$.

\end{lemma}

\demo{Proof:}
For $s\in S$ we have ${\cal H}_s =  \phi(\bigcup_{t\in \phi_T^{-1}(s)}{\cal G}_t)$ and  $\phi_T^{-1}({\cal H}_s) = \bigcup_{t\in \phi_T^{-1}(s)} {\cal G}_t$, hence
\[
 \phi(W) \cap {\cal H}_s = \bigcup_{t\in \phi_T^{-1}(s)} \phi({\cal G}_t \cap W).
\]
If $\sigma\in N\cap M$ and $t\in \phi_T^{-1}(s)$, then
$$
\phi({\cal G}_t \cap W)=\sigma\phi({\cal G}_t \cap W)=\phi(\sigma {\cal G}_{t} \cap \sigma W)=\phi({\cal G}_{\sigma t} \cap W).
$$
It follows that
$$
\bigcup_{t\in \phi_T^{-1}(s)} \phi({\cal G}_t \cap W)=\bigcup_{\sigma \in R} \phi({\cal G}_{\sigma t_0} \cap W), \hspace{.5cm}t_0\in \phi_T^{-1}(s),
$$
where $R$ is a (finite) system of representatives of $G/(N\cap M)$ in $G$. In particular,
$$
{\cal H}_s=\bigcup_{\sigma \in R' } \phi({\cal G}_{\sigma t_0})\,,
$$
$R'$ a system  of representatives of $G/N$ in $G$, hence we proved (i).

Since for any $t\in \phi_T^{-1}(s)$ the surjective homomorphism of pro-$\c$-groups ${\cal G}_t\surjr{} {\cal H}_s$ is closed, we see that $\phi(W) \cap {\cal H}_s$ is closed in ${\cal H}_s$. Furthermore, since  $\phi$ is strict, we have
\[
 \phi(W)\cap {\cal V} = \phi ({\cal U}\cap W)
\]
and  $\phi_{|{\cal U}}: {\cal U}\to {\cal V}$ is a  surjection. Since  ${\cal U}$ and ${\cal V}$ are  compact, $\phi_{|{\cal U}}$ is closed. It follows that $\phi(W)\cap {\cal V}$ is closed in ${\cal V}$.
 Thus  $\phi(W)$ is a closed subset of $({\cal H}, {\cal V})$.
\enddemo

\subsection{Corestricted bundles over a one-point compactification of a discrete set}

\indent

Let $T=T_0\ktcup\{*\}$ be the one-point compactification of a discrete set $T_0$.
Let $( G_t )_{t \in T_0 }$ be a  family of
pro-$\c$-groups indexed by the elements of the discrete set $T_0$. Consider the set
$$
{\cal G}:= \Tcup_{t \in T_0} G_t \, \, \ktcup \{ \ast \}\,,
$$
where $\{ \ast \}$ is the group with one element, equipped with the following topology: $G_t \subseteq {\cal G}$
(together with its profinite topology) is open in ${\cal G}$ for all $t$,
and for every open neighborhood $V \subseteq T$ of $\ast \in T$,
let\vspace{-.2cm}
$$
\Tcup_{t \in V} G_t\, \, \ktcup \{ \ast \}
$$

\smallskip\noindent
be an open neighborhood of $\ast\in {\cal G}$.
Then  $({\cal G},pr,T)$, where $pr$ is  the continuous map
$$
pr: {\cal G} \rightarrow T; \, \, G_t \ni g_t \mapsto t,\, \, \ast \mapsto \ast
$$
 is a compact bundle of pro-$\c$-groups.

\bigskip
 Now we consider the topology of a corestricted bundle over $T$. Let
 $$
 ({\cal G},{\cal U})=(\Tcup_{t \in T_0} {\cal G}_t\ktcup  \{*\},\Tcup_{t \in T_0} {\cal U}_t\ktcup \{*\})
 $$
  be a corestricted bundle of pro-$\c$-groups over $T$ with respect to the compact subbundle  ${\cal U}$.
 Then a subset $V$ of $({\cal G},{\cal U})$ is open if and only if the following holds:

\vspace{.2cm}
\Item{(i)} If $*\in V$, then ${\cal U}_t\subseteq V$ for almost all $t\in T_0$,
\Item{(ii)} $V\cap {\cal G}_t$ is open in ${\cal G}_t$ for all $t\in T_0$.

\vspace{.2cm}\noindent
Indeed, assume $V$ is open in $({\cal G},{\cal U})$, and so (ii) holds, and $W=V\cap {\cal U}$ is open in ${\cal U}$.  Let  $S=\{t\in  T| {\cal U}_t\subseteq W\}$. Then
$$
T\backslash S=\{t \in T\,|\, {\cal U}_t \nsubseteq W \}=pr_{{\cal U}}({\cal U}\backslash W).
$$
Since  ${\cal U}$ is compact, and  so ${\cal U}\backslash W$ is compact, it follows that $ T\backslash S$ is compact, hence $S$ is open in $T$.
If $*\in V$, then  $S$ contains $*$, thus $ T\backslash S$ is finite, i.e.\ (i) holds.
Conversely, assume that  (i) and (ii) holds. If $*\in V$,  then $ T\backslash S$ is finite. It follows that $W={\cal U}\backslash \ttcup_{ T\backslash S}({\cal U}_t\backslash W_t)$ is open in ${\cal U}$.  If $t\neq *$, then ${\cal U}_t$ is open in ${\cal U}$, and so every open subset of ${\cal U}_t$ is open in ${\cal U}$, thus $W=V\cap {\cal U}$ is open in ${\cal U}$, if $*\notin V$.

\bigskip

We will prove that under certain conditions a fibrewise surjective strict morphism $\phi: (\cal{G},\cal{U})\pfeil(\cal{H},\cal{V})$ of corestricted bundles is open, where $(\cal{H},\cal{V})$ is a corestricted bundle over the one-point compactification of a discrete set. We need
the following

\begin{lemma}\label{compactopenmapping}
 Let $(\Cal{U}, pr_{\Cal{U}}, T),\ (\Cal{V}, pr_{\Cal{V}}, S)$ be compact bundles of pro-$\c$-groups\linebreak where $T$ is a profinite space and $S=S_0\ktcup\{*\}$ is the one-point compactification of the discrete set $S_0$ such that $\Cal{V}_*=\{*_{\Cal{V}}\}$ is the group with one element. Assume that

  {
\smallskip
\Item{(i)}
$\phi: \Cal{U}\pfeil \Cal{V}$ is a fibrewise surjective morphism,

 \smallskip
\Item{(ii)}
 $\phi_T: T\pfeil S$ is an open map and

  \smallskip
\Item{(iii)}
 $\#\phi^{-1} (*_{\Cal{V}})=1$, i.e.\ $*_{\Cal{V}}$ has a unique preimage.

 }
 \smallskip\noindent
Then $\phi$ is open.
\end{lemma}

\demo{Proof:}
Let $W\subseteq \Cal{U}$ be an open subset, and so $W_t=W\cap \Cal{U}_t$ is open in $\Cal{U}_t$ for all $t\in T$. Since $\phi$ is fibrewise surjective, the homomorphisms $\phi_t: \Cal{U}_t\surjr{} \Cal{U}_{\phi_T(t)}$ are open. Hence it follows that for all $s\in S$ the set\vspace{-.2cm}
\[
 \phi(W)_s = \phi(W)\cap \Cal{V}_s = \bigcup_{t\in \phi_T^{-1}(s)} \phi_t(W_t)
\]
is open in $\Cal{V}_s$. In order to prove that $\phi(W)$ is open in $\Cal{V}$, by lemma (\ref{compfamily}) it remains to show that for a closed subset $X\subseteq \Cal{V}$, the set
\[
 S_X = \{s\in S\ |\ X_s\subseteq \phi(W)_s\}
\]
is open in $S$. Again by lemma (\ref{compfamily}), the set
\[
 T_X = \{t\in T\ |\ \phi^{-1}(X)_t\subseteq W_t\}
\]
is open in $T$. Hence, since $\phi_T$ is open, it follows that $\phi_T(T_X)$ is open in $S$. Obviously $\phi_T(T_X)\subseteq S_{X}$. Since $S=S_0\ktcup\{*\}$ is the one-point compactification of a discrete space, it remains to note that $*\in \phi_T(T_X)$ if $*\in S_X$. In fact, this is guaranteed by the assumption $\#\phi^{-1} (*_{\Cal{V}})=1$.
\enddemo

We are now able to prove the following

\begin{proposition}\label{openmappingtheorem}
 Let $(\Cal{G}, \Cal{U})$ be a corestricted bundle of pro-$\c$-groups over a profinite space $T$ and $(\Cal{H},\Cal{V})$ be a corestricted bundle over $S$ where $S=S_0\ktcup\{*\}$ is the one-point compactification of the discrete set $S_0$ such that $\Cal{H}_* = \{*_{\Cal{H}}\}$ is the group with one element. Let\vspace{-.2cm}
 $$
 \phi: (\Cal{G}, \Cal{U})\pfeil (\Cal{H},\Cal{V})
 $$
 be a fibrewise surjective strict morphism of corestricted bundles such that

   {
\smallskip
\Item{(i)} $\phi_T: T\pfeil S$ is an open map and

 \smallskip
\Item{(ii)}$\#\phi^{-1} (*_{\Cal{H}})=1$, i.e.\ $*_{\Cal{H}}$ has a unique preimage.

}
\smallskip\noindent
 Then $\phi$ is open.
\end{proposition}

\demo{Proof:}
Let $W\subseteq \Cal{G}$ be an open subset. For all $t\in T$, the set $W_t\subseteq \Cal{G}_t$ is open. Since $\phi$ is fibrewise surjective, the homomorphisms $\phi_t: \Cal{G}_t\surjr{} \Cal{H}_{\phi_T(t)}$ are open and hence for all $s\in S$ the set
\[
 \phi(W)_s = \phi(W)\cap \Cal{H}_s = \bigcup_{t\in \phi_T^{-1}(s)} \phi_t(W_t)
\]
is open in $\Cal{H}_s$. It remains to show that $\phi(W)\cap \Cal{V}$ is open in $\Cal{V}$. Note that since $\phi$ is strict and surjective, we have
\[
 \phi(W) \cap \Cal{V} = \phi(W\cap \Cal{U}).
\]
By lemma (\ref{compactopenmapping}), the induced map $\phi_{|_\Cal{U}}: \Cal{U}\surjr{} \Cal{V}$ is open which implies the claim.
\enddemo

One should remark that the assumption (i) in the above ``open mapping result'' is necessary by the following general observation: If $\phi\colon ({\cal G}, {\cal U})\to ({\cal H}, {\cal V})$ is an open morphism of corestricted bundles over $T$ and $S$, respectively, then $\phi_T\colon T\to S$ is open. In fact, let $T'\subseteq T$ be an open subset and let $W=\phi (\ttcup_{t\in T'} {\cal G}_t)\cap {\cal V}$ which is an open subset of ${\cal V}$ by assumption. Now lemma (\ref{compfamily}) implies that the set $\phi_T (T') = \{s\in S\,|\,\{ 1_t\}\subseteq W_t\}$, where $1_t$ denotes the unit element in ${\cal V}_t$, is open in S.
\bigskip

The following lemma shows that the canonical projections of projective limits of profinite spaces are open provided the transition maps are. More precisely, we have the following

\begin{lemma}\label{onepointcompactifiedlimit}
 Let $\{S_i,\rho_{ij}\}_I$ be a projective system of profinite spaces such that all transition maps $\rho_{ij}: S_i\to S_j, i\ge j$ are surjective. Let $S=\varprojlim_{j\in I} S_j$ and denote by $\rho_i$ the canonical surjection $S\twoheadrightarrow S_i,\ i\in I$. Then the following holds:

\smallskip
\Item{(i)} If all transition maps $\rho_{ij}$ are open, then the maps $\rho_i$ are open, too.

\smallskip
\Item{(ii)} If $S_i = S_{0i}\ktcup \{*_i\}$ is the one-point compactification of a discrete set $S_{0i}$ and $\rho_{ij}^{-1}(\{*_j\}) = \{*_i\}$ for all $i\ge j$, then the maps $\rho_{ij}$ and $\rho_i$ are open.
\end{lemma}

\demo{Proof:}
 Assertion (ii) follows from (i), since the maps $\rho_{ij}$ are obviously open. In order to prove (i), it is sufficient to show that the set $\rho_i(V\cap S)$ is open in $S_i$, where\vspace{-.2cm}
\[
 V=\prod_{j\in I} V_j\subseteq \prod_{j\in I} S_j
\]
with open and closed subsets $V_j\subseteq S_j$ such that $V_j=S_j$ for all $j\notin I_0$, where $I_0$ is a finite non-empty subset of $I$. For any $i\in I$ we define a subseteq $W_i\subseteq S_i$ as follows: Let $k\in I$ such that $k\ge i$ and $k\ge j$ for all $j\in I_0$ and set
\[
 W_i:=\rho_{ki}\big( \bigcap_{j\in I_0} \rho_{kj}^{-1} (V_j)\big) \subseteq V_i\subseteq S_i.
\]
Since $I$ is directed, it follows that $W_i$ is independent from the choice of $k$. Assuming that $\rho_{ij}$ is open for all $i\ge j$, $W_i$ is an open and closed subset of $S_i$. We claim that $W_i=\rho_i(V\cap S)$. It is easy to see that $\rho_{ij}$ maps $W_i$ surjectively onto $W_j$ for all $i\ge j$. Let $W=\varprojlim_{j\in I} W_j$ and denote by $\tilde{\rho_i}$ the canonical surjection $W\twoheadrightarrow W_i$. Since $W\subseteq V\cap S$, we have
\[
 W_i=\tilde{\rho}_i(W) \subseteq \rho_i(V\cap S).
\]
Conversely, let $x=(x_j)_I\in V\cap S$ and let $k\in I$ such that $k\ge i$ and $k\ge j$ for all $j\in I_0$. Then $\rho_{ki}(x_k)=x_i$ and $\rho_{kj}(x_k)=x_j\in V_j$ for all $j\in I_0$, i.e.\ $x_k\in \bigcap_{j\in I_0} \rho_{kj}^{-1} (V_j)$. Therefore $\rho_i(x)=x_i\in W_i$ which finishes the proof.
\enddemo

\begin{theorem} \label{projlim3}
Let $\{({\cal G}_i,{\cal U}_i),T_i,\phi_{i,j}\}_I$ be a projective system of corestricted bundles of pro-$\c$-groups, where $ T_i =T_{0i}\ktcup \{*_i\}$ is the one-point
compactification of a discrete set $T_{0i}$ and
 $({\cal G}_i)_{*_i} = \{*_{{\cal G}_i}\}$  for all ${i\in I}$. Let $T=\varprojlim_{i\in I}  T_i$ and
$$
({\cal G},{\cal U})=\varprojlim_{i\in I}({\cal G}_i,{\cal U}_i)
$$
as defined in section \ref{coreslimitsection}. Assume that the following holds.

{

\smallskip
\Item{(i)} The morphisms $\phi_{ij}$ are fibrewise surjective and strict.

\smallskip
\Item{(ii)} For all $i\ge j$ the transition maps satisfy $\phi_{ij}^{-1}(\{*_j\}) = \{*_i\}$.

}

\smallskip\noindent
Then
$$
\freeproduct_{T} {\cal G}=\varprojlim_{i\in I} \freeproduct_ {T_i} ({\cal G}_i,{\cal U}_i).
$$
\end{theorem}

\demo{Proof:} Using (\ref{openmappingtheorem}) and (\ref{onepointcompactifiedlimit}), we see that the canonical maps
$\phi_i\colon {\cal G}\twoheadrightarrow ({\cal G}_i,{\cal U}_i)$ are open. Now ({\ref{projlim1}}) (with $G=1$) and its following remark give the desired result.
\enddemo

In proposition (\ref{subbundles}) we have seen that every corestricted free product \linebreak ${\freeproduct}_{T}({\cal G},{\cal U})$ can be written as the completion of the unrestricted free product ${\freeproduct}_{T}({\cal G},{\it 1})$. In the case where $T$ is the one-point compactification of a discrete set $T_0$, we will show that ${\freeproduct}_{T}({\cal G},{\cal U})$ can also be obtained in a different way as projective limit involving unrestricted free products of quotient bundles with fibres ${\cal G}_t/{\cal U}_t$.
\bigskip

Assume that ${\cal U}_t$ is a normal subgroup of ${\cal G}_t$ for all $t\in T$. If $S$ is a finite subset of $T_0$, then  ${\cal U}^S$ denotes the compact subbundle $\ttcup_{t \in T\backslash S} {\cal U}_t$ of ${\cal U}$. According to (2.1), we have a canonical (open) surjection of bundles
$$
({\cal G},{\cal U})\twoheadrightarrow ({\cal G}/{\cal U}^S,{\cal U}/{\cal U}^S)=({\cal G}/{\cal U}^S,{\it 1}),
$$
where as sets
$$
{\cal G}/{\cal U}^S=\Tcup_{t \in S} {\cal G}_t\ktcup\Tcup_{t \in T\backslash S} {\cal G}_t/{\cal U}_t.
$$
It follows that we obtain an isomorphism of bundles
$$
({\cal G},{\cal U})\liso\varprojlim_S ({\cal G}/{\cal U}^S,{\it 1}),
$$
where $S$ runs trough the finite subsets of $T_0$. Thus we get the following proposition, see also \cite{N} Satz (2.2):

\begin{proposition} \label{projlim2}
Let $T=T_0\ktcup\{*\}$ be the one-point compactification of a discrete set $T_0$ and let $({\cal G},{\cal U})$ be a corestricted bundle of pro-$\c$-groups over $T$. Let ${\cal U}_t$ be  normal  in ${\cal G}_t$ for all $t\in T$. Then
$$
\freeproduct_{t\in T}({\cal G}_t,{\cal U}_t)=\varprojlim_S \big (\freeproduct_{t\in S}{\cal G}_t*\freeproduct_{t\in T\backslash S}({\cal G}_t/{\cal U}_t,{\it 1})\big ).
$$
\end{proposition}

\demo{Proof:}
If $S_1\subseteq S_2$ are finite subsets of $T_0$, then the morphism of bundles
$$
\phi_{S_2,S_1}: ({\cal G}/{\cal U}^{S_2},1)\twoheadrightarrow ({\cal G}/{\cal U}^{S_1},{\it 1})
$$
is fibrewise surjective and
 $\phi_{S}\colon ({\cal G},{\cal U})\twoheadrightarrow ({\cal G}/{\cal U}^S,{\it 1})$ is open. Using ({\ref{projlim1}}) (with $G=1$), we obtain the desired result.
\enddemo

\section{Corestricted free products of  families}

\indent

Let $G$ be a pro-$\c$-group and let $(G_t)_{t \in T}$ and  $(U_t)_{t \in T}$
be families of closed subgroups of $G$ indexed by the points of a profinite space
$T$, where $U_t$ is a closed subgroup of $G_t$ for every $t\in T$.  Assume that
$(U_t)_{t \in T}$ is a continuous family, i.e.\
$$
{\cal U} =\{(g,t)\in G \times T\,|\, g \in U_t \}
$$
is a compact bundle over $T$, see (\ref{family})(i).
Let
$$
{\cal G}  = \{(g,t)\in G \times T\,|\, g \in G_t \}
$$
and let $pr_{\cal G}$ be the restriction to ${\cal G}$ of the projection
$G \times T\twoheadrightarrow T$. Let  ${\cal G}$ be equipped with the following topology:
a set $V\subseteq {\cal G}$ is open if and only if

{
\medskip
\Item{(i)} $V\cap {\cal U}$ is open in ${\cal U}$,
\Item{(ii)}$V\cap {\Cal G}_t$ is  an open subset of $G_t$ for all $t\in T$.

}
\medskip\noindent
(We identify the fiber ${\cal G}_t=(G_t,t)$ with $G_t$.)  Observe that this topology on ${\cal G}$ is finer than the topology which is induced by the topology of the constant (compact) bundle $(G\times T, pr, T)$. We define the maps $m, \iota, e$ by
restricting  the corresponding maps from the constant bundle to ${\cal G}$.

\begin{lemma}\label{base} With the notation and assumptions as above let $x\in {\cal G}$ and let $V$ be an open neighborhood of $x$. Then there exists an open neighborhood $V_0$ of $x$ of the following form:
\Item{(i)} if $x\in U_t$, then $V_0=(xN\times S)\cap {\cal G} $, where $N$ is an open normal subgroup of $G$ and $S$ open in $T$,
and  $V_0\cap {\cal U}\subseteq V\cap {\cal U}$,
\Item{(ii)} if $x\in G_t\backslash U_t$, then
$V_0=xN_t$, where $N_t$ is an open normal subgroup of $G_{t}$, and $V_0\subseteq V$.
\end{lemma}

\demo{Proof:}
 If $x\in U_t$, then $x\in W=V\cap {\cal U}$. Since $W$ is open in the compact bundle  ${\cal U}$ (which equipped with the induced topology of the compact constant bundle $G\times T$), there exists an open normal subgroup $N$ of $G$ and an open subset $S$ of $T$ such that
 $(xN\times S)\cap {\cal U}\subseteq W$.
Then $V_0=(xN\times S)\cap {\cal G}$ is an open neighborhood of $x\in {\cal G}$ and $V_0\cap {\cal U}\subseteq V\cap {\cal U}$.

\medskip
  Assume that $x\in G_t\backslash U_t$ and let $N_t$ be an open normal subgroup of $G_t$ such that $xN_t\subseteq V_{t}=V\cap G_t$.
We may assume that $N_t$ is small enough such that $xN_t\cap U_{t}=\varnothing$. It follows that
the set $V_0=xN_t\subseteq G_{t}\subseteq {\cal G}$ is open in ${\cal G}$. Indeed, for  $t'\neq t$ the set $V_0\cap G_{t'}=\varnothing$ is open in $G_{t'}$,
the set $V_0\cap G_{t}=xN_t$ is open in $G_{t}$ and $V_0\cap {\cal U}=\varnothing$ is open in ${\cal U}$.
\enddemo

\begin{lemma}\label{subgroup}
Let $G$ be a pro-$\c$-group, $U_0\subseteq U$ are  closed subgroups of $G$ and $g$ an element of $G$. Let  $V$ be an open subset of $G$ containing $gU_0$ such that
$$
 gU_0= V\cap gU.
$$
 Then there exists an open subgroup $H$ of $G$ such that $gU_0\subseteq gH\subseteq V$ and $ U_0=H\cap U$.
\end{lemma}

\demo{Proof:} Since $U_0=\bigcap_{i\in I}H_i$, where $H_i$ is an open subgroup of $G$ for every $i\in I$, we have $G\backslash V\subseteq \bigcup_{i\in I}(G\backslash gH_i)$. Since $G\backslash V$ is closed in $G$, hence compact, and the sets $G\backslash gH_i$  are open ($H_i$ is open and closed), we get $G\backslash V\subseteq \bigcup_{i=1}^n(G\backslash gH_i)$, and so\vspace{-.2cm}
$$
gU_0 \subseteq \bigcap _{i=1}^n gH_i\subseteq V.
$$

\vspace{-.2cm}\noindent
Thus $H=\bigcap _{i=1}^n H_i$ has the desired properties.
\enddemo

\begin{theorem}\label{corfamily}
With the notation as above, the space  ${\cal G}=({\cal G},{\cal U})$  is  a corestricted bundle of pro-$\c$-groups over $T$
with respect to ${\cal U}$, and we have continuous inclusions
$$
 \xymatrix{{\cal U} \ar@{^(->}[r]&({\cal G},{\cal U}) \ar@{^(->}[r]^\phi &G\times T.
 }
$$
\end{theorem}

\demo{Proof:} Since the injective map $({\cal G},{\cal U})\injr{\phi} G\times T$ is continuous and  $G\times T$ is a totally disconnected Hausdorff space, the same is true for the space $({\cal G},{\cal U})$.

Obviously, the map $pr_{\cal G}\colon {\cal G}\pfeil T$ is continuous: if $S\subseteq T$ is open, then $pr_{\cal G}^{-1}(S)\cap {\Cal G}_t$ is empty or  equal to $G_t$ for all $t\in T$ and $pr_{\cal G}^{-1}(S)\cap {\Cal U}=pr_{\Cal U}^{-1}(S)$ is open in ${\Cal U}$.
In particular, it follows that a fiber $G_t$ is closed in $({\cal G},{\cal U})$
and that the topology on $G_t$ induced by the topology of ${\cal G}$ is the pro-$\c$ topology of $G_t$.

Now we prove that the multiplication  $m\colon {\cal G} \times_T {\cal G} \rightarrow {\cal G}$ is continuous.
Let $(a,b)\in {\cal G} \times_T {\cal G}$ and let
$t_0\in T$ such that $a,b\in G_{t_0}$. Let $V=\ttcup_{t\in T}V_t\subseteq {\cal G}$ be an open neighborhood of $m(a,b)=ab$.
We consider the following two cases.

1. Let $a\notin U_{t_0}$ or $b\notin U_{t_0}$. Assume that $a\notin U_{t_0}$ and let $N_{t_0}$ be an open normal subgroup of $G_{t_0}$ such that $abN_{t_0}\subseteq V_{t_0}$.
We may assume that $N_{t_0}$ is small enough such that $aN_{t_0}\cap U_{t_0}=\varnothing$. It follows that
the set $V_a:=aN_{t_0}\subseteq G_{t_0}\subseteq {\cal G}$ is an open neighborhood of $a$ in ${\cal G}$, see lemma (\ref{base}).

Let $V_b:=\ttcup_{t\neq t_0}G_t\ktcup (bN_{t_0})$. Then $V_b$ is open in $({\cal G},{\cal U})$, since $V_b\cap {\cal U}=\ttcup_{t\neq t_0}U_t\ktcup (bN_{t_0}\cap U_{t_0})$ is open in
${\cal U}$ and $V_b\cap G_t$ is open in $G_t$ for all $t\in T$. Furthermore, $m(V_a,V_b)=abN_{t_0}\subseteq V_{t_0}\subseteq V$.

2. Let $a, b\in U_{t_0}$.
Using lemma (\ref{base}), we replace $V$ by $(abN\times S)\cap V$, where $N$ is an open normal subgroup  of $G$ and  $S$ an open subset of $T$, i.e.\ we may assume that
$$
 W= V\cap {\cal U}=(abN\times S)\cap {\cal U}=\Tcup_{t\in S}(abN\cap U_t).
$$
Thus
$
V_t\cap U_t= abN\cap U_t
$
for every $t\in S$. Let
$$
 W_a=\Tcup_{t\in S}(aN\cap U_t)\hspace{.5cm}\hbox{and}\hspace{.5cm} W_b= \Tcup_{t\in S}(bN\cap U_t).
$$
Then $W_a$ and $W_b$ are open neighborhoods of $a$ and $b$ in ${\cal U}$, respectively, and
 $m_{\cal U}(W_a,W_b)\subseteq W$. Let
$$
S_a:=\{t\in S| aN\cap U_t\neq\varnothing\} \,\,\hbox{ and }\,\,S_b:=\{t\in S| bN\cap U_t\neq\varnothing\}.
$$
For every $t\in S_a$ resp.\ $t\in S_b$ there are elements $n_t\in N$  resp.\ $m_t\in N$ such that  $an_t\in U_t$ and $bm_t\in U_t$ and
$$
aN\cap U_t=an_t(N\cap U_t) \,\,\hbox{ resp.\ }\,\,bN\cap U_t=am_t(N\cap U_t).
$$
Let $t\in S_a\cap S_b$. Using lemma (\ref{subgroup}) and
$$
an_tbm_t(N\cap U_t)=abN\cap U_t=V_t\cap U_t=an_tbm_tU_t\cap V_t,
$$
there exists an open subgroup $H_t$ of $G_t$ such that
$$
an_tbm_t(N\cap U_t)\subseteq an_tbm_t H_t\subseteq V_t
$$
and $N\cap U_t= H_t\cap U_t$. Let
$$
\begin{array}{lcl}
V_a&=&\Tcup_{t\in S_a\cap S_b} an_t(H_t)^{bm_t}\,\,\ktcup\,\Tcup_{t\in S_a\backslash S_b} aN,\\
V_b&=&\Tcup_{t\in S_a\cap S_b} bm_tH_t\,\,\,\ktcup\,\,\Tcup_{t\in S_b\backslash S_a} bN,
\end{array}
$$
where $(H_t)^{bm_t}=(bm_t)H_t(bm_t)^{-1}$. Then $m(V_a,V_b)\subseteq V$. We show that the sets $V_a$ and $V_b$ are open in $({\cal G},{\cal U})$.
Obviously we have condition (ii). Furthermore,
$$
V_a\cap {\Cal U}=\Tcup_{t\in S_a}(aN\cap U_t)=W_a, \,\,\,\hbox{ and }\,\,\, V_b\cap {\Cal U}=\Tcup_{t\in S_b}(bN\cap U_t)=W_b,
$$
since
$$
 an_t(H_t)^{bm_t}\cap U_t= an_t(H_t\cap U_t)^{bm_t}= an_t(N\cap U_t)^{bm_t}=an_t(N\cap U_t)=aN\cap U_t.
$$

In order to prove that the inversion map
 $\iota\colon{\cal G} \rightarrow {\cal G}$ is continuous, let $a^{-1}\in  {\cal G}$ and  $V\subseteq {\cal G}$  an open neighborhood of $a^{-1}$. Then the set $V^{-1}=\{x\in  {\cal G}| x^{-1}\in V\}$ contains $a$, $\iota(V^{-1})=V$  and $V^{-1}$ is open in ${\cal G}$, since $V^{-1}\cap{\cal U}=\iota_{\cal U}^{-1}(V\cap{\cal U})$ is open in ${\cal U}$ as $\iota_{\cal U}$ is continuous, and $V^{-1}\cap {\Cal G}_t=(V\cap{\cal G}_t)^{-1}$ is  an open subset of $G_t$ for all $t\in T$.

\medskip
Finally,
 since the unit  $e$ is equal to the composition $T \lpfeil^{e_{\cal U}} {\cal U}\hookrightarrow {\cal G}$, the map $e$ is continuous.
\enddemo

The continuous map $\phi$ induces  a  homomorphism
$$
\phi_*: \freeproduct_T ({\cal G},{\cal U}) \lpfeil G.
$$

\begin{definition}\label{def-cor-freeprod2}
Let $(G_t)_{t\in T}$ and $(U_t)_{t\in T}$
be families of closed subgroups of a pro-$\c$-group $G$ indexed by the points of a profinite space
$T$.
Assume that $(U_t)_{t\in T}$ is a continuous family and  $U_t$ is a closed subgroup of $G_t$ for every $t\in T$.
Let $({\cal G},{\cal U})$
be the  corestricted bundle which is associated to these families.

We say that the  pro-$\c$-group $G$ is the  {\bf corestricted free pro-$\c$-product} of the family
$(G_t)_{t\in T}$ with respect to the continuous family $(U_t)_{t\in T}$ of closed subgroups of $G$
if $\phi_*$ is an isomorphism. We write\vspace{-.2cm}
$$
G=\freeproduct_ {t\in T} (G_t,U_t).
$$
\end{definition}

\noindent
{\bf Remarks:} We consider the situation given in definition (\ref{def-cor-freeprod2}).

\vspace{.2cm}\noindent
1. If   $({\cal G},{\cal U})$ is a corestricted bundle of pro-$\c$-groups over the one-point compactification $T$ is of a discrete set $T_0$. Then it can be considered as a corestricted bundle  associated to the families
$({\cal G}_t)_{t \in T }$ and $({\cal U}_t)_{t\in T}$
of closed subgroups of the  pro-$\c$-group $G={\freeproduct}_T({\cal G},{\cal U})$, see corollary (\ref{faser-injective1}); indeed, since the bundle ${\cal U}$ is compact, the family $({\cal U}_t)_{t\in T}$ is continuous, see (\ref{family})(ii).

\vspace{.2cm}\noindent
2. For all $t\in T$, let $U_t'\subseteq U_t$ be a closed subgroup such that $(U'_t)_{t\in T}$ is a continuous family. If $\Cal U\,'$ denotes the compact bundle associated to the family $(U'_t)_{t\in T}$, then we have a canonical surjection ${\freeproduct}_{T} ({\cal G}, {\Cal U\,'}) \surjr{} {\freeproduct}_{T} (\cal G, \cal U)$ and an isomorphism $ \varprojlim_{N} ({\freeproduct}_T ({\cal G},{\Cal U\,'}))/{N} \liso{}{\freeproduct}_T ({\cal G}, {\cal U})$, see proposition (\ref{subbundles}).

\vspace{.2cm}\noindent
3. If \ $t_0\in T$, then  the canonical map
$
\omega_ {G_{t_0}}:  G_{t_0}\pfeil\freeproduct_{t\in T} (G_t,U_t)
$
is an injective group homomorphism. This follows from the fact that the composition
$G_{t_0}\pfeil{\freeproduct}_{t\in T} (G_t,U_t)\pfeil G $ is injective.

\subsection{Abelianization of corestricted free products}

\indent

Let\vspace{-.5cm}
$$
\Resprod{t\in T_0} (A_t,B_t)=\big\{(a_t)_{t\in T_0}\in \sprod{t\in T_0} A_t\,\,|\,\, a_t\in B_t \hbox{ for almost all $t\in T_0$}\big\}
$$
be the restricted product over a discrete set $T_0$ of abelian  locally
compact groups $A_t$ with respect to  closed subgroups $B_t$.
The topology is given by the subgroups $V$ such that
\Item{(i)} $V\cap A_t$ is open in $A_t$ for all $t\in T_0$,
\Item{(ii)} $V\supseteq B_t$  for almost all $t\in T_0$.

\noindent
Then we call
$$
{\Resprod{t\in T_0}}^{\!\!c} (A_t,B_t):=\varprojlim_V \big (\Resprod{t\in T_0}(A_t,B_t)\big )/V,
$$
the {\it compactification} of $\resprod_{t\in T_0} (A_t,B_t)$, where $V$ runs through all open subgroups of finite index in $\resprod_{t\in T_0} (A_t,B_t)$. The canonical map
 $\resprod_{t\in T_0} (A_t,B_t)\pfeil\resprod_{t\in T_0}^c (A_t,B_t)$ has  dense image.

\medskip\noindent
We define the {\it discretization} of $\resprod_{t\in T_0} (A_t,B_t)$ by
$$
{\Resprod{t\in T_0}}^{\!\!d} (A_t,B_t):=
\varinjlim_W   W
$$
where $W$ runs through the finite subgroups of $\resprod_{t\in T_0} (A_t,B_t)$.
 If the subgroups $B_t$ of  $A_t$, $t\in T_0$, are open and compact, then $\resprod_{t\in T_0} (A_t,B_t)$ is locally compact. Using the equality
$
 (\resprod_{t\in T_0} (A_t,B_t))^\vee=\resprod_{t\in T_0} (A_t^\vee,(A_t/B_t)^\vee),
$
we obtain
$$
{\Resprod{t\in T_0}}^{\!\!d} (A_t,B_t)=(((\Resprod{t\in T_0} (A_t,B_t))^\vee)^c)^\vee \hspace{.5cm}\hbox{ and } \hspace{.5cm}
{\Resprod{t\in T_0}}^{\!\!c} (A_t,B_t)=(((\Resprod{t\in T_0} (A_t,B_t))^\vee)^d)^\vee,
$$
where $^\vee$ denotes the Pontryagin-dual.

\begin{proposition}\label{disc-comp} Let $T$ the one-point compactification of the discrete set $T_0$.

\noindent\smallskip
\Item{(i)} Let $A_t$,  $t\in T_0$,  be  discrete abelian torsion groups such that their exponents have a common finite bound. Then we have the equality of abstract groups
$$
\Resprod{t\in T_0}^{\!\!d} (A_t,B_t)=\Resprod{t\in T_0}(A_t,B_t).
$$
However, $\resprod_{t\in T_0}^d (A_t,B_t)$ is endowed with the discrete topology in contrast to $\resprod_{t\in T_0} (A_t,B_t)$.

\noindent\smallskip
\Item{(ii)} Let each $A_t$,  $t\in T_0$,  be a profinite abelian group. Then the canonical map $\resprod_{t\in T_0} (A_t,B_t)\hookrightarrow\resprod_{t\in T_0}^c (A_t,B_t)$ is injective. Setting $A_*=B_*=\{0\}$, then $(B_{t})_{t\in T}$ is a continuous familiy of closed subgroups of $\Resprod{t\in T_0}^{\!\!c} (A_t,B_t)$ and
$$
\freeproduct_ {t\in T} (A_t,B_t)\cong \Resprod{t\in T_0}^{\!\!c} (A_t,B_t)
$$
where ${\freeproduct}_ {t\in T} (A_t,B_t)$ is the corresponding corestricted free pro-abelian product.

\noindent\smallskip
\Item{(iii)}
Let ${\freeproduct}_{t\in T} (G_t,U_t)$
be  the  corestricted free pro-$\c$-product of the family
$\{G_t\}_{t \in T}$ with respect to the continuous family $\{ U_t\}_{t \in T}$ where $G_*=U_*=\{*\}$.
Then
$$
\big (\freeproduct_{t\in T} (G_t,U_t)\big )^{ab}={\Resprod{t\in T}}^{\!c}(G_t^{ab},\tilde U_t)=\varprojlim_N \big (\Resprod{t\in T}(G_t^{ab},\bar U_t)\big )/N,
$$
where $G_t^{ab}=G_t/[G_t,G_t]$ and   $\bar U_t=U_t[G_t,G_t]/[G_t,G_t]$
and $N$ runs through the subgroups of  $\resprod_{t\in T}(G_t^{ab},\bar U_t)$ of finite index such that $N\cap G_t^{ab}$ is open in $G_t^{ab}$ for all $t\in T$ and $N\supseteq \bar U_t$ for almost all $t\in T$.
\end{proposition}

\demo{Proof:} Under the assumption of (i) every element of $\resprod_{t\in T_0} (A_t,B_t)$ generates a finite subgroup. Thus (i) is obvious.

In order to prove (ii), we fix $t_0\in T_0$. If $V_{t_0}$ is an open subgroup of $A_{t_0}$, then $\tilde V_{t_0}:=V_{t_0}\times \resprod_{t\neq t_0} (A_t,B_t)$ is an open subgroup of $\resprod_{t\in T_0} (A_t,B_t)$ of finite index. If $V_{t_0}$ runs through a basis of neighborhoods of the unit of $A_{t_0}$, then the $t_0$-component of the intersection of the open subgroups $\tilde V_{t_0}$ is $\{0\}$. Varying $t_0$, it follows that  the intersection
$\bigcap V$ of all open subgroups $V$ of finite index in $\resprod_{t\in T_0} (A_t,B_t)$ is zero. This proves the first assertion of (ii). The second assertion of (ii) follows immediately from the definitions. Finally, (iii) follows from (ii) and (\ref{abelian}). In fact, noting that $\{s\}\subseteq T$ is open and closed in $T$ for all $s\in T_0$, by (\ref{zerlegung}) (ii) there exists a continuous splitting of $G_s\pfeil \freeproduct_{t\in T} (G_t,U_t)$.
\enddemo

\section{Cohomology of corestricted free products}

\indent

Now we consider the cohomology of a corestricted free pro-$\c$-product
$$
G=\freeproduct_{t\in T} (G_t,U_t)
$$
of a family  $(G_t)_{t\in T}$ with respect to a continuous family $(U_t)_{t\in T}$ of closed subgroups of a pro-$\c$-group $\tilde G$
in the  case where $T=T_0\ktcup \{*\}$ is the one-point
compactification of a discrete set $T_0$. By $\tilde U_t$ we denote  the normal closure of $U_t$ in $G_t$.

\begin{lemma}\label{normal-closure}
With the above notation $(\tilde{U}_t)_{t\in T}$ is a continuous family. Furthermore, if ${\cal G}=\ttcup_{t \in T_0} G_t\ktcup \{*\},\ {\cal U}=\ttcup_{t \in T_0} U_t\ktcup \{*\}$ and $\tilde{\cal U}=\ttcup_{t \in T_0} \tilde{U}_t\ktcup \{*\}$, then the (continuous) morphism of bundles $id: (\cal G, U)\to (G,\tilde{U})$ induces an isomorphism of pro-$\c$-groups
\[
 \freeproduct_{t\in T} (G_t, U_t) \liso \freeproduct_{t\in T} (G_t, \tilde{U}_t).
\]
\end{lemma}
\demo{Proof:}
Let $V$ be an open neighborhood of the identity in $\tilde{G}$. Then $V$ contains a normal subgroup $N$ of $\tilde{G}$. Since $\{U_t\}_{t\in T}$ is a continuous family,
the set
\[
 T(N)=\{t\in T\ |\ U_t\subseteq N\} = \{t\in T\ |\ \tilde{U}_t\subseteq N\}
\]
is an open subset of $T$ containing $\{*\}$. Since $T$ is the one-point compactification of the discrete set $T_0$, it follows that also the set
\[
 T(V)=\{t\in T\ |\ \tilde{U}_t\subseteq V\}
\]
containing $T(N)$ is open in $T$. Hence, the family $(\tilde{U}_t)_{t\in T}$ is also continuous and we have the continuous (not necessarily open) morphism of bundles $id: (\cal G, U)\to (G,\tilde{U})$. In order to show that the corresponding corestricted free products coincide, it remains to show that any morphism $({\cal G, U})\to H$, where $H$ is a finite $\c$-group, is continuous with respect to the topology of $({\Cal G},\tilde{\Cal U})$. However, this is a direct consequence of remark 2 following definition (\ref{def-cor-freeprod2}).
\enddemo

Let $A$ be  a discrete ${\freeproduct}_ {t\in T} (G_t,U_t)$-module. Let
$H^i_{nr}(G_t,A)$ be defined as the image of the inflation map $H^i(G_t/\tilde U_t,A^{\tilde U_t})\!\pfeil \!H^i(G_t,A)$.
It is easy to see that the map
$$
\res_H: H^i(G, A) \pfeil \prod_{t\in T}H^i(G_t, A)
$$
has image in ${\Resprod{T}}^d \!(H^1(G_t, A),H^1_{nr}(G_t,A))$, see  \cite{N}  \S 4.

\begin{theorem}\label{coh-freeprod} With the notation and assumptions as above let $A$ be a finite  ${\freeproduct}_ {t\in T} (G_t,U_t)$-module.
Then there is an exact sequence
$$
0\pfeil\! A/A^G \pfeil\!
{\Resprod{T}}^d \! (A/A^{G_t},A^{\tilde U_t}/A^{G_t})\pfeil
H^1(G, A) \pfeil\!
{\Resprod{T}}^d \!(H^1(G_t, A),H^1_{nr}(G_t,A))\!
\pfeil 0
$$
and an isomorphism\vspace{-.2cm}
$$
H^2(G, A) \liso
{\Resprod{T}}^d (H^2(G_t, A),H^2_{nr}(G_t,A)).
$$

\vspace{-.2cm}\noindent
If the cohomological dimension of $G_t/\tilde U_t$ is equal or less than $1$ for all $t\in T$, then
$$
H^i(G, A) \liso
{\bigoplus_{T}} H^i(G_t, A),\quad i\geq 3.
$$
\end{theorem}

\demo{Proof:} By lemma (\ref{normal-closure}), we may assume without loss of generality that $U_t=\tilde{U}_t$, i.e.\ $U_t$ is a normal subgroup of $G_t$ for any $t\in T$. If $i\leq 2$, the assertion   follows along the same lines as in  \cite{N} Satz (4.1); but one has to be careful with the topology and has to replace the restricted product by its discretization, see (\ref{disc-comp})(i).

If $i\geq 3$, we use dimension shifting: Let $A'$ be defined by the exact sequence $0\pfeil A\pfeil \Coind_G A\pfeil A'\pfeil 0$, where $\Coind_G A$ denotes the coinduced $G$-module consisting of all continuous functions from $G$ to $A$. By assumption we have $H^i_{nr}(G_t,A)=0$ for all $t\in T$ and $i\geq 2$. Thus the assertion follows since we get compatible isomorphisms $\phi\colon H^{i-1}(G, A')\kiso H^i(G, A)$ and $\phi_t\colon H^{i-1}(G_t, A')\kiso H^i(G_t, A)$, see also \cite{N} Satz (4.2).
\enddemo

Now we consider in the following case which has an application in number theory.
Let $I$ be a directed set and let\vspace{-.2cm}
$$
G=\varprojlim_{\lambda\in I}  G_\lambda,
$$

\vspace{-.2cm}\noindent
be a pro-$\c$-group given as projective limit of pro-$\c$-groups $G_\lambda$. Let
$
T=\varprojlim_{\lambda\in I}  T_\lambda,
$
where
$ T_\lambda =T_{0\lambda}\ktcup \{*_\lambda\}$ is the one-point
compactification of a discrete set $T_{0\lambda}$. Let $\{{\cal G},{\cal U}\}=
\{({\cal G}_\lambda,{\cal U}_\lambda)\}_\lambda$ be a projective system of corestricted bundles $({\cal G}_\lambda,{\cal U}_\lambda)$ over $T_\lambda$ with transition maps $\phi_{\mu\lambda}$, where
${\cal G}_\lambda=\ttcup_{t_\lambda \in T_{0\lambda}} G_{t_\lambda}\ktcup \{*_\lambda\}$ and  ${\cal U}_\lambda=\ttcup_{t_\lambda \in T_{0\lambda}} U_{t_\lambda}\ktcup \{*_\lambda\}$ and
$
U_{t_\lambda}\subseteq G_{t_\lambda}\subseteq G_\lambda
$
are closed subgroups, i.e.\ the diagrams
$$
\xymatrix{
U_{t_\mu}\ar@{^(->}[r]\ar[d]& G_{t_\mu}\ar@{^(->}[r]\ar[d]& G_\mu\ar[d]\\
U_{t_\lambda}\ar@{^(->}[r]& G_{t_\lambda}\ar@{^(->}[r]&G_\lambda
}
$$
commute for $\mu\geq \lambda$ and $\phi_{\mu\lambda}(t_\mu)=t_\lambda$.
Let
${\cal U}=\varprojlim_\lambda {\cal U}_\lambda$ and ${\cal G}=\varprojlim_\lambda {\cal G}_\lambda$ be the projective limits
 and let
$G_t=\varprojlim_\lambda G_{t_\lambda}$, $U_t =\varprojlim_\lambda U_{t_\lambda}$.
Assume that the following holds:

{

\smallskip
\Item{(i)} the morphisms $\phi_{\mu\lambda}$ are fibrewise surjective and strict,

\smallskip
\Item{(ii)} for all $\mu\ge \lambda$ the transition maps satisfy $\phi_{\mu\lambda}^{-1}(\{*_\lambda\}) = \{*_\mu\}$.

}

\smallskip\noindent
By theorem (\ref{projlim3}) we have
$$
\freeproduct_{ T} ({\cal G},{\cal U})
=\freeproduct_ {t\in T} (G_t,U_t)\liso\varprojlim_\lambda \freeproduct_ {t_\lambda\in T_\lambda} (G_{t_\lambda},U_{t_\lambda})
$$
and we have a canonical homomorphism
$$
\varphi: \freeproduct_ {t\in T} (G_t,U_t)\lpfeil G.
$$

\noindent\medskip
{\bf Remark:} In the number theoretical situation we have in mind, one can also argue with  theorem (\ref{corestricted-proj-lim}) in order to establish the isomorphism above.

\begin{proposition}\label{freeprod-limit} With the notation and assumptions as above let $p$ be a prime number and let $\c$ be the class of finite $p$-groups. Then the following are equivalent.

\medskip
\Item{(i)} $\freeproduct_ {t\in T} (G_t,U_t)\lriso{\varphi} G$ is an isomorphism.

\medskip
\Item{(ii)} The induced maps
 $$
 \varphi_*: H^1(G, \Z/p\Z) \kiso\varinjlim_{\lambda} {\Resprod{T_\lambda}}^d (H^1(G_{t_\lambda}, \Z/p\Z),H^1_{nr}(G_{t_\lambda},\Z/p\Z))
 $$
 $$
 \varphi_*: H^2(G, \Z/p\Z) \hookrightarrow\varinjlim_{\lambda} {\Resprod{T_\lambda}}^d (H^2(G_{t_\lambda}, \Z/p\Z),H^2_{nr}(G_{t_\lambda},\Z/p\Z))
$$
are bijective resp.\ injective.
\end{proposition}

\demo{Proof:} Using (\ref{coh-freeprod}), we have
$$\renewcommand{\arraystretch}{1.5}
\begin{array}{lcl}
H^i(\freeproduct_ {t\in T} (G_t,U_t),\Z/p\Z) &=&\ds\varinjlim_{\lambda}H^i(\freeproduct_ {t_\lambda\in T_\lambda} (G_{t_\lambda},U_{t_\lambda}),\Z/p\Z)\\
&=& \ds\varinjlim_{\lambda}{\Resprod{T_\lambda}}^d (H^i(G_{t_\lambda}, \Z/p\Z),H^i_{nr}(G_{t_\lambda},\Z/p\Z))
\end{array}
$$
for $i=1,2$. Now the usual argument, see \cite{NSW} (1.6.15), gives the result.
\enddemo


We have the following application in number theory. Let $p$ be a prime number and let $k$ be number field and $T$ a set of primes of $k$. We use the following notation.
$$
\begin{array}{ll}
k(p)& \hbox{is  the maximal $p$-extension of $k$,}\\
k^T& \hbox{is the maximal $p$-extension of $k$ which is completely decomposed} \\
&\hbox{at every prime of $T$,}\\
G(k(p)|k^T)& \hbox{is  the Galois group of the extension $k(p)|k^T$,}\\
G_\P(k)& \hbox{is the decomposition group of $G(k(p)|k)$ with respect to  $\P$,}\\
I_\P(k)& \hbox{is and inertia group of $G(k(p)|k)$ with respect to  $\P$,}
\end{array}
$$
where  $\P$ is an extension of a prime $\p$  of $k$  to $k(p)$.

\begin{theorem}\label{numberth}Let  $p$ be an odd prime
number and let $k$ be a number field of CM-type containing the group $\mu_p$ of all $p$-th roots of unity,  with maximal totally real subfield $k^+$, i.e.\
$k=k^+(\mu_p)$ is totally imaginary  and $[k:k^+]=2$.
Let
$$
T=\{\p\,|\,\p\cap k^+\hbox{ is inert in $k|k^+$}\}\cup \{\p|p\}.
$$
 Then the Galois group $G(k(p)|k^T)$ is the  corestricted free pro-$p$-product of the family
$(G_\P(k))_{\P\in {\cal T}}$ with respect to the continuous family $(I_\P(k))_{\P\in {\cal T}}$, i.e.\ the canonical map
$$
\freeproduct_ {\P\in {\cal T}} (G_\P(k),I_\P(k))\liso G(k(p)|k^T)
$$
is an isomorphism. Here ${\cal T}=\varprojlim_K \bar T_K$ is the projective limit of the one-point compactifications $\bar T_K$ of the discrete sets $T_K$ of all prolongations  of $T$ to $K$ and $K$ runs through all finite Galois extensions inside $k^T|k$.
\end{theorem}

This follows from (\ref{freeprod-limit}) and the results of \cite{W}: (2.4), (2.5), (2.2).

\medskip
\begin{tabbing}
\hspace{1cm}\= Mathematisches Institut\hspace{.8cm}  e-mail: gaertner \,\,(at) mathi.uni-heidelberg.de \\
\>der Universit\" at Heidelberg  \hspace{1.7cm} wingberg (at) mathi.uni-heidelberg.de\\
\>Im Neuenheimer Feld 288  \\
\>69120 Heidelberg  \\
\>Germany
\end{tabbing}


\begin{thebibliography}{999}

\bibitem{GR} Gildenhuys, D., Ribes, L.\ {\it A Kurosh subgroup theorem for free pro-$\c$-products of pro-$\c$-groups}. Trans.\ Amer.\ Math.\ Soc. {\bf 186} (1973)
309--329
\bibitem{Har} Haran, D.\ {\it On closed subgroups of free products of profinite groups}.
Proc.\ London Math.\ Soc.\ (3) {\bf 55} (1987) 266--298
\bibitem{Mel} Melnikov, O.\ V.\, {\it Subgroups and homology of free products of profinite groups} (in Russian).
Izv.\ Akad.\ Nauk SSSR {\bf 53} (1989). English translation
in Math.\ USSR Izv.\ {\bf 34} (1990) 97--119
\bibitem{N}  Neukirch, J.\, {\it Freie Produkte pro-endlicher Gruppen und ihre
Kohomologie}. Archiv der Math.\ {\bf 22} (1971) 337--357
\bibitem{NSW} Neukirch, J., Schmidt, A., Wingberg, K. \,{\it Cohomology of
 Number Fields, 2nd edition.} Springer 2008
 \bibitem{P} Pontryagin, L.\ S.\, {\it Topological Groups}.  New York, London, Paris
1966
\bibitem{W} Wingberg, K.\ {\it Sets of Completely Decomposed Primes in Extensions of  Number Fields}. Preprint 2013

\end{thebibliography}
\end{document}